\documentclass[12pt]{article}
\usepackage{amssymb,amsmath, amsthm}
\usepackage{graphicx,epsfig}
\usepackage{caption}
\usepackage{subcaption}
\usepackage{enumitem}
\usepackage{epstopdf} 
\usepackage{xcolor}
\usepackage{float}
\usepackage{rotating}
\usepackage{tikz}
\usepackage{dsfont}
\usepackage{comment} 
\usepackage{tabularx}
\usepackage{hyperref}
\usepackage{upgreek}
\usepackage[margin=3cm]{geometry}
\usepackage[english]{babel}
\usepackage[square,numbers]{natbib} 
\usepackage{booktabs} % For professional-looking tables
\usepackage{algorithm}
\usepackage{algpseudocode}
\hypersetup{
    colorlinks=true,
    linkcolor=cyan,
    filecolor=cyan,      
    urlcolor=cyan,
    citecolor = cyan
    }

 \numberwithin{equation}{section}

\newtheorem{assumption}{A.}
\newtheorem{assumptionB}{B.}

\begin{document}
\newcounter{eqnarray}%[section]
\newtheorem{t1}{Theorem}[section]
\newtheorem{d1}{Definition}[section]
\newtheorem{c1}{Corollary}[section]
\newtheorem{l1}{Lemma}[section] \newtheorem{r1}{Remark}[section]
\newtheorem{e1}{Example}[section]
\newtheorem{p1}{Proposition}[section]
\newtheorem{a1}{Result}[section]

\title{A novel characterization of structures in smooth regression curves: from a viewpoint of persistent homology   }

\author{\small 
	Satish Kumar \\
	\small IIT Kanpur\\
	\small Department 
 of Mathematics and Statistics \\
	\small  Kanpur 208016, India\\
	{\small Email: satsh@iitk.ac.in }\\
	\and
	\small Subhra Sankar Dhar \\
	\small  IIT Kanpur\\
	\small   Department of Mathematics and Statistics \\
	\small Kanpur 208106, India\\
	{\small Email: subhra@iitk.ac.in}\\
}
\maketitle

\begin{center}
    \textbf{Abstract} 
\end{center} 
We characterize structures such as monotonicity, convexity, and modality in smooth regression curves using persistent homology. Persistent homology is a key tool in topological data analysis that detects higher-dimensional topological features such as connected components and holes (cycles or loops) in the data. In other words, persistent homology is a multiscale version of homology that characterizes sets based on the connected components and holes. We use super-level sets of functions to extract geometric features via persistent homology. In particular, we explore structures in regression curves via the persistent homology of super-level sets of a function, where the function of interest is --- the first derivative of the regression function. In the course of this study, we extend an existing procedure of estimating the persistent homology for the first derivative of a regression function and establish its consistency. Moreover, as an application of the proposed methodology, we demonstrate that the persistent homology of the derivative of a function can reveal hidden structures in the function that are not visible from the persistent homology of the function itself. In particular, we characterize structures such as monotonicity, convexity, and modality, and propose a measure of statistical significance to infer these structures in practice. Finally, we conduct an empirical study to implement the proposed methodology on simulated and real data sets and compare the derived results with an existing methodology. 

\noindent\textbf{Keywords:} Topological data analysis, Persistent homology, Non-parametric regression, Regression curves, Level sets, Convexity, Monotonicity, Modality.        

\section{Introduction} \label{Introduction} 
The shape of the regression function plays a crucial role in postulating the hypothesis about the regression model. One of the ways of characterizing the shape of a function is by using derivatives of the function. Derivatives of a function provide important information about the shape of the function, such as monotonicity, convexity, critical points, inflection points, etc. Therefore, to gain insights into the shape of the regression function from the observed data, one can estimate the derivatives and use them to infer geometric features of the underlying true regression function. One such way of estimating the derivative is by estimating the regression function from suitable smoothing techniques, for example, the Nadaraya-Watson estimator, and then taking its derivative. There are various ways of estimating derivatives of a regression function (see, e.g., \cite{GasserSacnd}, \cite{DerLPR}, \cite{OptimalDer}). However, we use the derivative of the Nadaraya-Watson estimator as an estimate of the derivative of the regression function.   

There have been various attempts to explore the shape of the regression function from a geometric perspective; for instance, \cite{ShapeofReg} proposed a rank correlation coefficient to compare the shape of regression functions. In a different line of approach, motivated by the scale-space view from computer vision literature, \cite{SiZer} proposed a graphical device called SiZer to explore structures, such as peaks and valleys, in regression curves. SiZer analysis has been used in comparing two or more regression curves (see, e.g., \cite{WOOPARK20083954}, \cite{JCGS2014}). \cite{Dette2021} proposed a graphical device, based on the derivative and inverse of the derivative of the regression function, to compare the two regression functions up to a shift in the horizontal and/or vertical axis.  \cite{Dharwu} compared the time-varying quantile regression curves up to a horizontal shift.     

However, unlike the aforesaid works, we explore the shape of the regression function from a topological perspective. In particular, we attempt to describe the shape of the regression function in terms of connected components and holes, using the persistent homology (see, e.g., \cite{Edelsbrunner(2008)}, \cite{Ghrist(2007)}). One of the advantages of using topological features to describe the shape of a function is that these features provide qualitative global geometric information about the data and are stable under deformations. Hence, statistical inference based on these topological features, such as connected components and holes, is robust to various deformations in the data (see, e.g., \cite{Carlsson(2009)}, \cite{Carlsson(2014)}, \cite{Carlsson_Vejdemo-Johansson_2021}).

Now, let $f$ be a real-valued function. Then for a given $t \in \mathbb{R}$, the super-level set of $f$ at the level $t$ is defined as 
\begin{equation} \label{Dt}
    \displaystyle{ D_{t} \triangleq \{x \in \mathbb{R}: f(x) \geq t}\}.
\end{equation} 
The set $D_{t}$ contains important structural information about the function $f$, and therefore, is a fundamental object to analyze in many statistical procedures. For instance,  \cite{hartigan1975clustering} defines the number of clusters in a statistical population as the number of connected components of $D_{t}$, where the underlying function is the probability density function. \cite{MüllerSawitzki} proposed a test of multimodality of the probability density function based on the super-level sets. There is a vast literature on level-set estimation (see, e.g.,  \cite{Cuevas2000}, \cite{Tsybakov}), and numerous applications of super-level sets in statistical procedures (see, e.g., \cite{Cuevas1990}, \cite{Devroye1980}, \cite{Cuevas1997}).     

We analyze super-level sets of a function from a topological perspective. In particular, we take the function to be analyzed as the first derivative of the regression function and then calculate the persistent homology of its super-level sets. Since the regression function is usually not known in practice, we estimate the persistent homology using an estimate of the derivative of the regression function. In this work, we take the derivative of the Nadaraya-Watson estimator as an estimator of the derivative of the regression function. The validity of such estimators of persistent homology is due to the well-known stability theorem (see \cite{Cohen-Steiner2007}). The stability theorem provides an upper bound of the bottleneck distance between true persistent homology and the estimated persistent homology in terms of the sup-norm distance between the underlying function and the corresponding estimator. Thus, if the estimator of the underlying function is uniformly convergent, then the corresponding persistent homology is consistent with respect to the bottleneck metric. 

However, in the present context, we do not have any estimator of the derivative of the regression function, which is uniformly convergent. To circumvent this issue, we adopt the procedure proposed by \cite{Tdaconsistency} to estimate the persistent homology. This procedure is applied in a particular case where the estimator of the underlying functions is based on kernel estimation and the underlying functions are taken to be either density functions or regression functions.  \cite{Tdaconsistency} established the consistency of the proposed procedure up to the precision of $5 \epsilon$, where $\epsilon$ is any positive number involved in controlling topological features. The main advantage of this procedure is that it does not require uniform consistency of the estimator of the underlying function for the stability of the estimator of persistent homology.      
  
\subsection{Contribution}

The persistent homology of a function does not reveal local geometric features of its own accord. It only reveals global topological features such as connected components and holes in higher dimensions. It is worth mentioning here that the use of persistent homology to study the geometric structures of 1-manifolds is mainly limited to clustering structures. The use of persistent homology to infer local geometric structures is not so straightforward, and such procedures have been less developed in the literature. In this work, we study the local geometric features such as monotonicity, convexity, and modality of smooth regression functions not by studying their persistent homology but by exploring the persistent homology of their derivative. The methodology adopted here appears analogous to the procedure proposed by \cite{Tdaconsistency} to estimate the persistent homology of the derivative of the regression function, though they are entirely different too. In the present work, the goal is not only to recover the persistent homology of the true function but also to infer local geometric features of the underlying function using persistent homology. In contrast, the goal in \cite{Tdaconsistency} is to devise a new procedure to infer the persistent homology of density and regression functions. \cite{Tdaconsistency} established the consistency of the proposed estimator of persistent homology for density and regression functions using their kernel estimators. In a complementary fashion, we extend those results to the first derivative of the regression function and establish the consistency of the estimator of persistent homology up to precision $5 \epsilon$.

As an application of the proposed procedure, we demonstrate that it can be used to study the shape of regression functions. In particular, we characterize structures such as monotonicity, convexity, and modality in smooth regression curves and propose a significance measure of these structures observed in the data. Moreover, one can adopt the proposed procedure to compare two or more regression curves based on the structures such as monotonicity, convexity, and modality. In particular, one can perform cluster analysis for regression curves based on these structures. One of the real-life applications of this work is that, given the regression curves of population growth of different subcontinents, one can investigate the shape of the first derivative of regression functions to compare the population growth patterns over time. 

\subsection{Challenges}

Most of the techniques in the literature first obtain an estimator of the underlying function converging uniformly to its population counterpart and then use this estimator to estimate the persistent homology. As mentioned previously, such procedures are supported by the stability results for persistent homology (see Theorem \ref{stability theorem}). 

There are various methods for estimating the derivative of the regression function, but we are not aware of any result that establishes uniform consistency of the estimator. Therefore, most of the existing methods of estimating the persistent homology are ruled out in this case. To fix this, we use the procedure proposed by \cite{Tdaconsistency} to estimate the persistent homology of the derivative of the regression function. We restrict our analysis only to the first derivative of the regression function due to the complex form of the estimator, which makes calculation intractable. Even with the first derivative, due to the complex form of the estimator, it is not feasible to use Bernstein’s inequality while calculating the false positive and false negative error of the set $\{ X_{i}: \hat{f}_{n}(X_{i}) \geq t, i = 1, \ldots n \}$ with respect to $D_{t}$ (see, Equation \ref{Dt}), where $X_{i}$'s are observed data, and $\hat{f}_{n}$ is a kernel estimator of the underlying function $f$. Therefore, we use Hoeffding's inequality, which provides less sharp bounds than Bernstein’s inequality used in \cite{Tdaconsistency}. 

Moreover, in applying the proposed methodology in practice, a proper choice of kernel bandwidth is important. In particular, to assess the statistical significance of observed bars computed from the proposed methodology, the choice of kernel bandwidth depends on the unknown quantities, depending on the underlying true derivative. To address this issue, we propose an algorithm to estimate the true derivative values, which in turn allows estimation of the unknown quantities needed for an appropriate kernel bandwidth selection.

\subsection{Outline}

The rest of the article is organized as follows. Section \ref{Preliminaries} gives the necessary topological background to elaborate the notion of persistent homology. Section \ref{Problem Formulation and Main Results} formulates the exact problem and the main results of the article. Section \ref{Quantitative evaluation} presents a measure of statistical significance of observed bars.   Section \ref{Simulation studies} presents characterizations of the local structures, such as monotonicity, convexity, and modality, and provides some illustrative examples to implement the proposed methodology. Section \ref{Real data study} contains illustrations of the proposed methodology on two real data sets. Section \ref{Conclusion} consists of some concluding remarks. Finally, the proof of the main theoretical result is supplied in the appendix in Section \ref{Appendix} along with some other numerical and theoretical results.    

\section{Preliminaries} \label{Preliminaries}
Persistent homology (see, e.g., \cite{Edelsbrunner(2008)}, \cite{Ghrist(2007)}, \cite{ZomorodianCarlsson}) is a prominent tool in topological data analysis (TDA) that provides a multi-scale summary of topological features present in the data. Persistent homology is a multiscale extension of homology, which is one of the refined topological invariants from algebraic topology (see, e.g., \cite{Munkres(1984)}, \cite{Hatcher(2002)}). Heuristically, homology characterizes sets based on the number of connected components and holes using some algebraic machinery. In what follows, we are concerned with the homology of simplicial complexes referred to as simplicial homology. Simplicial complexes are combinatorial representations of topological spaces made up of points, lines, edges, triangles, and their higher-dimensional analogs. In the following subsections, we briefly describe the necessary background and definitions related to algebraic topology.           

\subsection{Simplicial homology} \label{Simplicial homology} 
In this subsection, we first define simplicial complexes, and then we describe the procedure to compute the homology of simplicial complexes. 

\begin{d1} [\textbf{Simplex}]
For any integer $k \geq 0$,  a $k$-simplex is defined to be the convex hull of $k + 1$ affinely independent points $ \mathbb{V} = \{x_{0}, x_{1}, \ldots, x_{k}\}$ in $\mathbb{R}^{d}, d \geq 1$. The elements of $\mathbb{V}$ are called vertices of the $k$-simplex. 
\end{d1}

In particular, points are 0-simplices, edges are 1-simplices, triangles are 2-simplices, and tetrahedra are 3-simplices. 

\begin{d1} [\textbf{Faces of a simplex}]
    Given a $k$-simplex and its vertex set $\mathbb{V}$, the faces of the $k$-simplex are defined to be the simplices spanned by the subsets of $\mathbb{V}$.
\end{d1}

\begin{d1} [\textbf{Geometric simplicial complex}]
A geometric simplicial complex or a simplicial complex $\mathcal{K}$ in $\mathbb{R}^d$ is a collection of simplices satisfying the following conditions:

\begin{itemize}
    \item If $\tau$ is any face of a simplex $\sigma$ in $\mathcal{K}$, then $\tau \in \mathcal{K}$.

    \item If $\tau$ and $\sigma$ are two simplices in $\mathcal{K}$, then either $\tau \cap \sigma = \emptyset$ or $\tau \cap \sigma$ is common face of both $\tau$ and $\sigma$.  
\end{itemize}
    
\end{d1}

Note that, when the vertex set of a geometric simplicial complex $\mathcal{K}$ is known, one can describe $\mathcal{K}$ completely and combinatorially by the list of its simplices. This leads to the notion of an abstract simplicial complex. 

\begin{d1} [\textbf{Abstract simplicial complex}]
Let $\mathcal{X} = \{ x_{0},\ldots, x_{n} \}$ be a finite set, and $\Delta$ be a collection of non-empty subsets in $\mathcal{X}$. Then $\Delta$ is said to be an abstract simplicial complex if it satisfies the following properties: 
\begin{enumerate}
    \item $\{x_{i} \} \in \Delta$ for all $i = 0, \ldots, n$.

    \item If $\tau\in\Delta$ and $\sigma\subseteq\tau$ such that $\sigma \neq \varnothing$, then $\sigma\in\Delta$.
\end{enumerate}
    
\end{d1} 

The set $\mathcal{X}$ is also called the vertex set of $\Delta$ and the elements of $\mathcal{X}$ are referred to as vertices. For an integer $k \geq 0$, any subset $\sigma \in \Delta$ with $k+1$ vertices is said to be a k-simplex or a k-dimensional simplex. Any subset $\tau $ of a k-simplex $\sigma$ with k vertices is called a face of $\sigma$. The dimension of $\Delta$ denoted as dim($\Delta$), is defined to be the maximum dimension of simplices in $\Delta$. That is, $ \text{dim}(\Delta) \triangleq \text{ max } \{ \text{dim}(\sigma): \sigma \in \Delta\}$. 

Now, we define the following algebraic structure on simplices in $\Delta$ to define holes or cycles in different dimensions of a topological space. First, we need to assign orientations; an ordering on the vertices of simplices in $\Delta$. Let $ \sigma = \{ x_{i_{0}}, \ldots, x_{i_{k}}\}$ be a k-simplex in $\Delta$, where $i_{0},\ldots,i_{k}$ is a permutation of $0, 1, \ldots, k$. We say two orderings of its vertex set are equivalent if they differ from one another by an even permutation. Thus, orderings of the vertices of $\sigma$ fall into two equivalent classes. Each of these classes is referred to as an orientation of $\sigma$. A 0-simplex has only one orientation. For any $k \geq 0$, we denote an oriented simplex as $[\sigma] = [x_{i_{0}}, \ldots, x_{i_{k}}$]. 

We denote the set of all k-simplices in $\Delta$ as $\Delta_{k}$. A k-chain on $\Delta$ is defined to be a formal linear combination of k-simplices in $\Delta_{k}$. Let $n_{k} = |\Delta_{k}|$, and $\gamma$ denotes a k-chain on $\Delta$, that is $\gamma \triangleq \sum\limits_{i=1}^{n_{k}} \alpha_{i}\sigma_{i}$, where the coefficients $\alpha_{i}$ belong to some field $\mathbb{F}$. The space of k-chains is a vector space over the field $\mathbb{F}$, and defined as $$C_{k}(\mathcal{X},\mathbb{F}) \triangleq \left\{\gamma = \sum\limits_{i=1}^{n_{k}} \alpha_{i}\sigma_{i}: \alpha_{i}\in\mathbb{F}\right\}.$$

Now, we need to map the simplices in $\Delta_{k}$ into a basis of $C_{k}(\mathcal{X},\mathbb{F})$ such that the map preserves orientations of simplices in $\Delta_{k}$. Let $\Delta^{\pi}_{k}$ denote the set containing all the simplices in $\Delta_{k}$ in all possible orientations. Any oriented simplex $[\sigma]$ in $\Delta^{\pi}_{k}$ is denoted as $[x_{\pi(0)}, \ldots, x_{\pi(k)}]$, where $\pi = \{\pi(0), \ldots, \pi(k)\}$ is a permutation of $0, \ldots, k$. Then, for any simplex $\sigma_{i} \in \Delta_{k}$, define the map $T_{k} : \Delta^{\pi}_{k} \xrightarrow{} C_{k}(\mathcal{X},\mathbb{F})$ as follows: $$ T_{k}([x^{i}_{\pi(0)}, \ldots, x^{i}_{\pi(k)}]) = (-1)^{P(\pi)} \boldsymbol{e}_{i}, $$ where $P(\pi)$ is the parity of the permutation $\pi$ and $\boldsymbol{e}_{i}$ is the standard basis in $\mathbb{R}^{n_{k}}$.   

A linear transformation $\partial _{k} : C_{k}(\mathcal{X},\mathbb{F}) \xrightarrow{} C_{k-1}(\mathcal{X},\mathbb{F}) $ is called the boundary operator. The $k^{th}$ boundary operator is a $n_{k-1} \times n_{k} $ matrix $\partial _{k}$ whose $i^{th}$ column $(\partial_{k})_{i}$ is defined as follows:
$$ (\partial_{k})_{i} = \sum \limits_{\sigma \in \Delta_{k-1}} T_{k-1} (\sigma),$$ where $\sigma$ is a face of $\sigma_{i}$. The subspaces $ B_{k}(\Delta)\triangleq \text{Im } \partial_{k +1} \text{ and } Z_{k}(\Delta)\triangleq \text{ ker }\partial_{k}$ are the set of k-boundaries and the set of k-cycles, respectively. The $k^{th}$ homology of $\Delta$ is defined to be the quotient vector space given by 
$$ H_{k}(\Delta) \triangleq \text{ ker }\partial_{k}/ \text{Im } \partial_{k +1}.$$ 
 
One of the refined topological invariants that are derived from the homology groups is the Betti numbers, and these Betti numbers provide a quantitative description of topological features in a space.  
 For any integer $k \geq 0$, we define the $k^{th}$ Betti number, $\beta_{k}$, of $\Delta$ as $\beta_{k}\triangleq \text{dim }(H_{k}(\Delta))$. Betti numbers are important topological invariants that count the number of connected components and holes in a topological space. In particular, $0^{th}$ Betti numbers count the number of connected components, and $1^{st}$ order Betti numbers count the number of one-dimensional loops or cycles. 

In practice, homology computation of a topological space $\mathcal{X}$ from the finite set of observations $\mathcal{X}_{n} = \{ x_{1},\ldots, x_{n} \}$ from it, involves converting the data $\mathcal{X}_{n}$ into simplicial complexes. Simplicial complexes are constructed by identifying simplices of different dimensions, which requires defining the proximity of data points at some threshold, say $r \geq 0$. The proximity of data points at a threshold $r$ is defined using some metric (say, Euclidean metric when data is embedded into an Euclidean space). The two frequently used constructions of simplicial complexes are the $\check{C}$ech complex and the Vietoris–Rips complex (see, e.g., \cite{Ghrist(2007)}).

The choice of proximity threshold $r$ is unknown in practice, and homology at a particular threshold is noisy. Therefore, one computes topological features (connected components and holes) at multiple thresholds and tracks the evolution of these features as thresholds vary. Then one considers a topological feature to be significant if it persists for a longer duration. This leads to the idea of persistent homology. 

\subsection{Persistent homology} \label{Persistent homology}
Let $\mathcal{X}$ be a space, then a filtration of $\mathcal{X}$ is a nested family of spaces $\{\mathcal{X}_{t} : t \in \mathbb{R} \}$ such that $ \mathcal{X} = \bigcup_{t \in \mathbb{R} } \mathcal{X}_{t}$. There are two types of filtrations often used in TDA. First, filtrations built on top of a point cloud, for example, the families of Rips or $\check{C}$ech complex indexed over a set of non-negative real numbers. Second, filtrations obtained through level sets of a function $f: \mathcal{X} \xrightarrow{} \mathbb{R}$, for example, one can take $\mathcal{X}_{t} = \{ x \in \mathcal{X}: f(x) \geq t\}$ to form a filtration of $\mathcal{X}$, such a filtration is referred to as a super-level set filtration of $f$. In practice, the index parameter $t$ is called the scale parameter. The homology of the spaces $\mathcal{X}_{t}$ varies with the scale parameter t. For example, as t increases, new connected components may be born, or existing connected components may disappear. Similarly, new holes or cycles are formed or filled up as t increases. Thus, as one moves through the filtration, new features are born, and existing features die. The birth time of a topological feature corresponds to the value of the scale parameter where it appears for the first time, and the death time of a feature corresponds to the value of the scale parameter where it dies or merges with another existing feature. Persistent homology tracks this process and records the pair of birth and death times of each topological feature in a multiset of points in $\mathbb{R}^2$, called persistence diagrams. One can refer to  \cite{Edelsbrunner(2008)}, \cite{Edelsbrunner(2010)}, \cite{ZomorodianCarlsson}, \cite{Carlsson(2009)} for rigorous treatment of persistent homology.    

Given a filtration $\{\mathcal{X}_{t}\}_{t \in \mathbb{R}}$ of a topological space $\mathcal{X}$, let, for any integer $k \geq 0$, $PH_k(\mathcal{X})$ denote the persistent homology of $\mathcal{X}$ that tracks the evolution of k-dimensional topological features as the scale parameter varies. Further, let $dgm_{k}$ denote the $k^{th}$ persistence diagram; a multiset that consists of pairs of the birth and death time of k-dimensional topological features in $PH_k(\mathcal{X})$, and let (b, d) correspond to the birth and death pair of a feature in $PH_k(\mathcal{X})$, then the persistence of the feature is defined as d - b. Features that have higher persistence are generally considered significant features, and those that correspond to lower persistence are regarded as topological noise. Another equivalent and useful representation of features is via persistence barcodes (see \cite{barcode}, \cite{Ghrist(2007)}). The persistence barcode is a set of intervals, where endpoints correspond to the birth and death times of a feature. A bar in the persistence barcode is the length of the interval [b, d]. Individual bars corresponding to each feature are stacked on top of each other so that the resulting structure resembles a barcode. 

\begin{figure}[!htt]
    \centering
    \includegraphics[width = \linewidth]{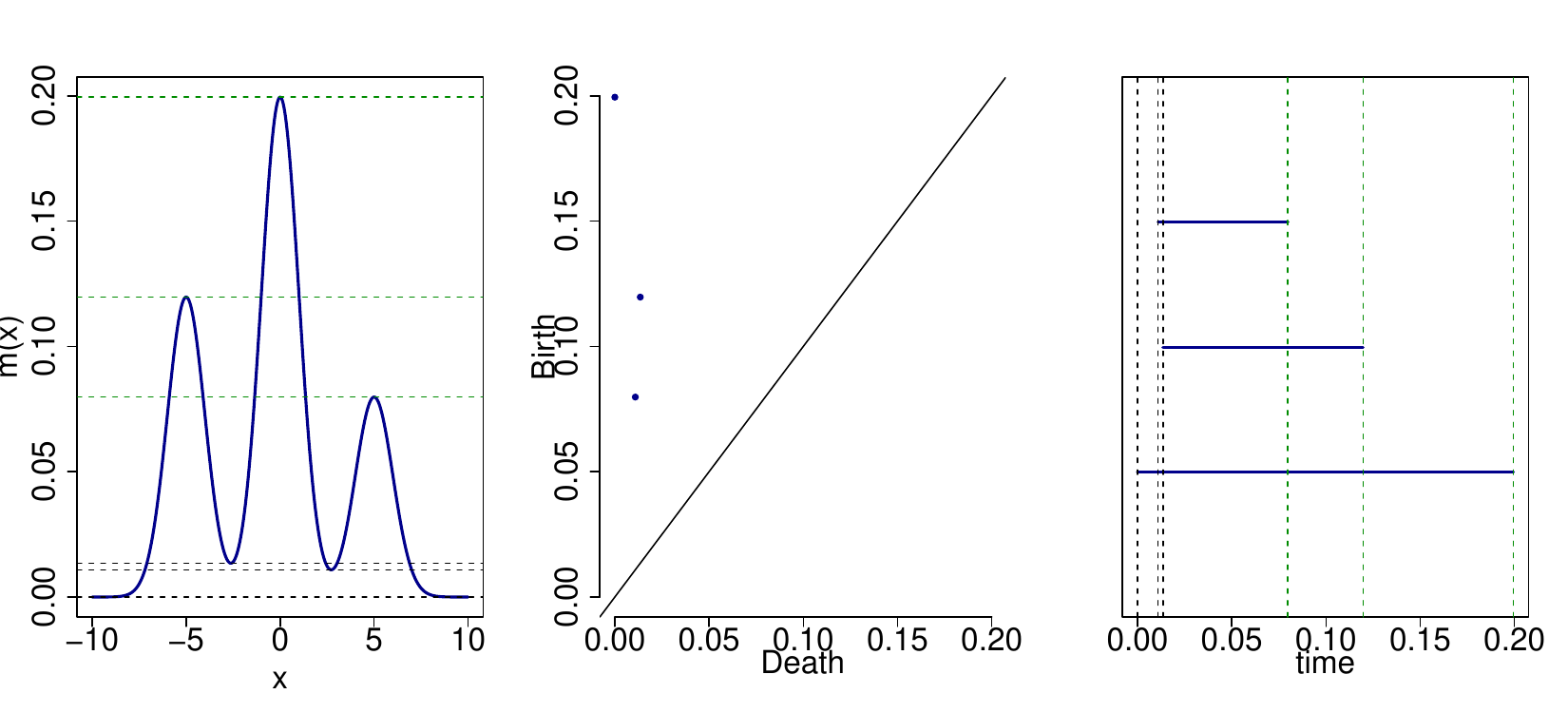}
    \caption{The graph of the function (left) $m(x) =0.5 p_1(x)+ 0.2 p_2(x)+0.3 p_3(x)$, where $p_1(.), p_2(.)$ and $p_3(.)$ denote the probability density functions of normal distributions, each with variance 1 and means 0, 0.5, and -0.5,  respectively. The persistence diagram of the super-level set filtration of $m$ is in the middle, and the persistence barcode of $m$ is on the right. The right endpoints of the bars in the barcode correspond to the local maxima, and the left endpoints correspond to the local minima of $m$. Green dotted lines indicate local maxima, and the black dotted lines indicate local minima. }
    \label{fig: persistent homology} 
\end{figure}

Figure \ref{fig: persistent homology} depicts the persistent homology of the super-level set of a function via the persistence diagram and the persistence barcode. The leftmost figure is the graph of a non-negative function with three modes. We consider the super-level set filtration of the function to study the structures via persistent homology. Note that at the level $t > 0.20$, the super-level set of the function is empty, hence there are no non-trivial connected components or cycles to study in this case. For $t \leq 0$, the super-level set of the function equals the entire domain [-10, 10], implying that there is only one connected component in the super-level set. The interesting and non-trivial topological features occur when $t \in (0, 0.20] $. Note that at $t = 0.20$, the super-level set is $\{0\}$, indicating the birth of a connected component. As we decrease $t$, new elements are added to the super-level set, resulting in the birth and death of topological features. The green dotted lines indicate the levels at which a new connected component is born, while the black dotted lines indicate the levels at which a connected component gets merged with the existing connected components. This process of evolution of features is summarized in the persistence diagram and the persistence barcode. Now, let $[d, b]$ denote an interval in the persistence barcode; then $b$ corresponds to the birth of a topological feature, and $d$ corresponds to the death of a topological feature for a super-level set filtration. For the super-level set filtration, the birth of a connected component corresponds to a local maximum, and the death of a connected component corresponds to a local minimum of the function.              

In this work, we are mainly concerned with the persistent homology of the super-level set filtrations of a real-valued function $f$, which is unknown. Therefore, one observes finitely many values of $f$ and estimates the persistent homology using an estimate $\hat{f}$ of $f$ (see, e.g., \cite{Fasy2014}, \cite{Chazal2011}, \cite{Chazal2013}, \cite{Bubenik2010}, \cite{Chung2009}, \cite{phillips_et_al2015}). The consistency of these procedures can be established using the well-known stability theorem for persistence diagrams (see \cite{Cohen-Steiner2007}). 

We need to define the following terminology before giving the formal statement of the stability theorem.

\begin{d1} [\textbf{Morse function}] \label{Morse function} 
  Let $\mathcal{M}$ be a d-dimensional manifold, and $f: \mathcal{M} \xrightarrow{} \mathbb{R}$ be a twice differentiable function. A point $p \in \mathcal{M} $ is said to be a critical point of $f$ if $\nabla f(p) = 0$, and the value $f(p)$ is called a critical value of $f$. A critical point $p$ is said to be non-degenerate if the Hessian matrix evaluated at $p$ is non-singular. Then the function $f$ is called a Morse function if all its critical points are non-degenerate, and all its critical values are distinct.
\end{d1}  

\begin{d1} [\textbf{Tame function}] \label{Tame functions}
 
Let $\mathcal{M}$ be a d-dimensional manifold, $f : \mathcal{M} \xrightarrow{} \mathbb{R}$ be a continuous function, and for any $t \in \mathbb{R}$, define a super level set $D_{t} = \{x \in \mathcal{M} : f(x) \geq t \}$. Then, we say that $t$ is a homological regular value if there exists $\epsilon > 0$ such that for every $ a \leq b $ in $( t - \epsilon, t + \epsilon)$, the map $H_{k}(D_{b}) \xrightarrow{} H_{k}(D_{a}) $, induced by the inclusion of super level sets, is an isomorphism for every $k \geq 0$. Otherwise, we say that $t$ is a homological critical value. Then the function $f$ is said to be tame if it has finitely many homological critical values, and $dim(H_{k}(D_{{t}})) < \infty$ for all $t$ and $k$. Note that any Morse function on a compact manifold is tame.
\end{d1} 

\begin{d1} [\textbf{Triangulable space}]
 Given a geometric simplicial complex $\mathcal{K}$, the underlying space of $\mathcal{K}$ is the union of the simplices of $\mathcal{K}$. The underlying space of $\mathcal{K}$, when viewed as a subset of an Euclidean space, is called a polytope. We say a topological space $X$ is triangulable if there exists a finite simplicial complex with its polytope homeomorphic to $X$.   
\end{d1}

\begin{d1} [\textbf{Bottleneck distance}] \label{Bottleneck distance}  Let, for $k\geq 0$, $PH_{k}(f)$ and $PH_{k}(g)$ denote the $k$-th persistent homology of the super-level set filtrations constructed via the continuous real-valued functions $f$, and $g$, respectively. Further, let $dgm_{k}(f)$ and $dgm_{k}(g)$ denote the $k$-th persistence diagrams corresponding to $PH_{k}(f)$ and $PH_{k}(g)$, respectively. Then, the bottleneck distance between $PH_{k}(f)$ and $PH_{k}(g)$ is defined as:
$$\delta_{B} ( PH_{k}(f), PH_{k}(g) ) = \inf_{\gamma \in \Gamma } \sup_{p \in dgm_{k}(f)  } || p - \gamma(p)||_{\infty},$$
where the set $\Gamma$ consists of all the bijections $\gamma: dgm_{k}(f) \cup \text{diag} \xrightarrow{} dgm_{k}(g) \cup \text{diag} $. Here,  diag = $\{(b, d) : b = d\}$, and $||.||_{\infty}$ denotes the sup-norm in $\mathbb{R}^2$.
\end{d1}
Note that the cardinality of both persistence diagrams may differ. Therefore, to define bijections between two persistence diagrams, we add the diagonal line (diag) to each persistence diagram, where each point has infinite multiplicity. The features with longer lifetimes (persistence) are far from the diagonal line, while features with shorter lifetimes lie near the diagonal line. The bottleneck distance between two persistence diagrams can be considered a measure of the maximum amount of displacement required to make two diagrams identical. Now, we provide the formal statement of the stability theorem.  
\begin{t1} [\textbf{Stability theorem} \cite{Cohen-Steiner2007}] \label{stability theorem} 
    Let $f$ and $g$ be real-valued continuous tame functions defined over the same triangulable space, and for any integer $k \geq 0$, $PH_{k}(f)$ and $PH_{k}(g)$ denote the corresponding $k^{th}$ persistent homology of $f$ and $g$, respectively. Then we have,
    $$ \delta_{B}\left( PH_{k}(f), PH_{k}(g)\right) \leq || f - g ||_{\infty} .$$
\end{t1}

\section{Problem formulation and main results} \label{Problem Formulation and Main Results}
Let $(X_{1},Y_{1}), \ldots ,(X_{n},Y_{n})$ be $n$ independent replications of a bivariate random vector $(X, Y)$. Further, let F denote the joint distribution of $(X, Y)$ with the joint probability density function $f$, and $p$ denote the marginal density of $X$. Consider the following non-parametric regression model:
\begin{equation} \label{Non parametric model}
     Y_{i} = m(x_{i}) + u_{i}, i = 1, \ldots , n,
\end{equation}
where $x_{i}$'s are observed values of the random variables $X_{i}$'s, $u_{i}$'s are unobserved errors assumed to be independent and identically distributed random variables with $\mathbb{E}(u_{i}) = 0$ and $var(u_{i}) = \sigma^2 (< \infty) $, and $m$ is the nonparametric regression function defined as $\displaystyle{m(x) \triangleq \mathbb{E}(Y|X = x )}$.    

 We aim to infer the geometric properties of the regression function $m$ using its first-order derivative. Let $m_{1}$ denote the first order derivative of the regression function $m$, we wish to summarize geometric features of $m$ through the super-level sets of $m_{1}$ which are  defined as, 
\begin{equation} \label{Level set}
    D_{t} = \{ x \in \mathbb{R} : m_{1}(x) \geq t \}, t \in \mathbb{R}.
\end{equation}
Note that $\{D_{t}\}_{t \in \mathbb{R}}$ is a nested family of super-level sets as $t$ decreases from $\infty$ to $-\infty$ that is for any $t_{1} < t_{2}$, we have $D_{t_{2}} \subset D_{t_{1}}$. We use persistent homology for the filtration $\{D_{t}\}_{t \in \mathbb{R}}$ to summarize the geometric features of $m_{1}$. Also note that $m_{1}$ is unknown since $m$ is unknown, and hence, we estimate persistent homology of $m_{1}$ from the observed data $\mathcal{X}_{n} \triangleq \{(X_{1},Y_{1}), \ldots (X_{n},Y_{n})\}$. The common procedure is to use an estimate $\hat{f}$ to estimate the persistent homology of an unknown function $f$. 

However, as mentioned previously, due to the uncertainty of stability results, we resort to the procedure developed in \cite{Tdaconsistency}. The main results in \cite{Tdaconsistency} present a robust estimator of homology for the level sets of density and regression functions. Then, using this robust estimator of homology, the procedure for computing persistent homology is extended. We need the following definition before defining the procedure formally.   
\begin{d1} \label{epsilon regular}
 Let $t \in \mathbb{R}$ be a level of $D_{t}$ (see, Equation \ref{Level set}), and $\epsilon \in (0, |t|/2)$. We say that $t$ is $\epsilon$-regular if the following holds:
\begin{align*}
\partial D_{t + 2\epsilon} \cap \partial D_{t + (3/2)\epsilon} &= \partial D_{t + (1/2)\epsilon} \cap \partial D_{t} = \partial D_{t} \cap \partial D_{t -(1/2)\epsilon}\\
& = \partial D_{t - (3/2)\epsilon} \cap \partial D_{t - 2\epsilon } = \varnothing,   
\end{align*}
where $\partial A$ denotes the boundary of the set A and $\varnothing$ denotes the empty set. In particular, if the underlying function $f$ is continuous on $f^{-1} [t - 2\epsilon, t + 2\epsilon]$, then $t$ is $\epsilon$-regular.   
\end{d1} 
  
The procedure proposed by \cite{Tdaconsistency} is summarized as follows:
\begin{enumerate}
    \item Let $f: \mathbb{R}^d \xrightarrow{} \mathbb{R} $ be a continuous tame (see, Definition \ref{Tame functions}) and bounded function. In this article $f = m_{1}$ and $d = 1$.  Consider a super-level set $D_{t}$ of $f$, and assume that the level $t$ is $\epsilon$-regular (see, Definition \ref{epsilon regular}). 

    \item Using the data $(X_{1}, Y_{1}), \ldots, (X_{n}, Y_{n})$ construct an estimator $\hat{f}_{n}$ using a kernel function $K$ with the bandwidth $h \equiv h(n)$, and define 
    \begin{equation} \label{Filtered set}
        \mathcal{X}^{t}_{n} \triangleq \{X_{i} : \hat{f}_{n}(X_{i}) \geq t, i = 1, \ldots, n  \}. 
    \end{equation}

    \item Now use $\mathcal{X}^{t}_{n}$ to construct an estimator of the super-level set $D_{t}$ as 
    \begin{equation} \label{hatDnh}
        \hat{D}_{t} (n, h) \triangleq \bigcup_{x \in \mathcal{X}^{t}_{n}} \mathbb{B}_{h}(x),
    \end{equation}
     where $h\equiv h(n)$ is the kernel bandwidth, and $\mathbb{B}_{h}(x)$ denotes a ball of radius $h$ centered around the point $x$. 

    \item Note that for any $\epsilon \in (0, |t|)$, we have $\hat{D}_{t + \epsilon} (n, h) \subset \hat{D}_{t - \epsilon} (n, h) $. This gives an inclusion map from $\hat{D}_{t + \epsilon} (n, h)\text{ to } \hat{D}_{t - \epsilon} (n, h)$.  

    \item The inclusion map $i : \hat{D}_{t + \epsilon} (n, h) \hookrightarrow \hat{D}_{t - \epsilon} (n, h)  $ induces a map in homology $$ i_{k} : H_{k}(\hat{D}_{t + \epsilon}(n, h)) \xrightarrow{} H_{k}(\hat{D}_{t - \epsilon}(n, h)) ,$$ where $k$ is any arbitrary degree of homology, and $h \equiv h(n)$ is the kernel bandwidth. 

    \item Now, use the image of the map $i_{k}$ as the robust estimator of homology denoted as $\hat{H}_{k}(t, \epsilon; n)$. Thus, for any arbitrary degree of homology $k$, we have    
    \begin{equation} \label{Robust estimator}
        \hat{H}_{k}(t, \epsilon; n) \triangleq \text{Im}(i_{k}).
    \end{equation}

    \item Define the following: 
    \begin{equation} \label{ Theorem constants}
        N_{\epsilon} \triangleq \sup_{x \in \mathbb{R}^d} \left\lceil f(x)/ 2 \epsilon \right\rceil, t_{max} = 2\epsilon N_{\epsilon} \text{ and } t_{i} = t_{max} - 2i\epsilon, i\in \mathbb{Z},      
    \end{equation}
    where $\left\lceil x \right\rceil$ denotes the smallest integer larger than or equal to $x$.

    \item Now, estimate the homology using $\hat{H}_{k}(t_{i}, \epsilon; n)$ for every $i \in \mathbb{Z}$, and compute the persistent homology for the following discrete filtration $\hat{\mathcal{D}}^{\epsilon} \triangleq \{\hat{D}_{t_{i}}(n, h)\}_{i\in \mathbb{Z}}$. This estimator of persistent homology is denoted as $\widehat{PH}^{\epsilon}_{k}(f)$.  
\end{enumerate}

 In this article, we establish the consistency of $\widehat{PH}^{\epsilon}_{k}(f)$ when $f = m_{1}$. We estimate $m_{1}$ using the Nadaraya-Watson estimator (see \cite{silverman1986density}) of the regression function $m$. That is, first, we estimate $m$ as
\begin{equation} \label{NW estimator}
    \hat{m}_{n}(x) = \frac{\sum\limits_{i = 1}^{n} K_{h}(x - X_{i}) Y_{i}}{\sum\limits_{i = 1}^{n} K_{h}(x - X_{i})},  
\end{equation}
where $K$ is a kernel function, $h \equiv h(n)$ is kernel bandwidth and $K_{h}(x) = K(x/h)$. In the following, we define an appropriate class of kernel functions and the required assumptions on the kernel functions.  

\begin{assumption} \label{A 1}
    The class of kernel functions is defined as:
$$\left\{K(x) = \frac{g(x)}{g(0)} : g(0) > 0, g \text{ is a  p.d.f } \text{ and }g(-x) = g(x) \in [0, g(0)], \forall x \in S\right\},$$
where p.d.f is an abbreviation of probability density function, and $S \subseteq [-1, 1]$ denotes the support of the density $g$.
\end{assumption}

\begin{assumption} \label{A 2}
    The kernel function $K$ is bounded away from 0, that is $\inf_{x \in S}K(x) = \delta$ where $\delta > 0$.
\end{assumption}

Now, we define the estimator of $m_{1}(x)$ as the first order derivative of $\hat{m}_{n}(x)$, that is 
\begin{equation} \label{Derivative of m}
    \hat{m}_{1, n}(x) = \frac{\mathrm{d}}{\mathrm{d}x}\hat{m}_{n}(x) = \frac{ \sum\limits_{i=1}^n Y_{i}K_{1,h}\left(x-X_{i}\right)}{\sum\limits_{i=1}^n K_{h}(x-X_{i})} - \frac{ \sum\limits_{i=1}^n Y_{i}K_{h}\left(x-X_{i}\right)\sum\limits_{i=1}^n K_{1,h}\left(x-X_{i}\right) }{\left(\sum\limits_{i=1}^n K_{h}(x-X_{i})\right)^2},
\end{equation}  
where $\displaystyle{K_{1,h}(x-X_{i}) = \frac{\mathrm{d}}{\mathrm{d}x} K_{h}(x-X_{i})}$ assuming that the first order derivative of the kernel function $K$ exists. In addition, we impose the following assumption on the kernel functions satisfying Assumptions A.~\ref{A 1} and A.~\ref{A 2}.

\begin{assumption} \label{A 3}
    For any real number $b > 0$, the derivative of kernel functions takes the form 

$$ \frac{\mathrm{d}}{\mathrm{d}x} K\left(\frac{x}{b}\right) = \frac{x}{b^2} \psi\left(\frac{x}{b}\right)K\left(\frac{x}{b}\right), $$ 
where $\psi$ is a function such that $ \tau_1\leq \psi(x) \leq \tau_2$, for all $x \in \mathbb{R}$. Here, $\tau_1, \tau_2 \in \mathbb{R}$ such that $\tau_1 \leq \tau_2$.

\end{assumption}
\begin{r1}
    The standard Gaussian and Cauchy kernels satisfy Assumption A.~\ref{A 3} with $\psi(x) = -1$, and $\psi(x) = 2\left(  1 + x^2\right)^{-1}$,  respectively, for $x \in \mathbb{R}$.    
\end{r1}

We need the following assumptions on the response variables, regressors, the regression function $m$, and its derivative $m_{1}$ to establish the consistency of $\widehat{PH}^{\epsilon}_{k}(m_{1})$.

\begin{assumptionB} \label{B 1}
    The marginal density of $X$ denoted as $p$ is tame (see, Definition \ref{Tame functions}) and bounded. Moreover, $X_{min}< X < X_{max}$ almost surely, where $X_{min}$ and $X_{max}$ are real numbers.
\end{assumptionB} 

\begin{assumptionB} \label{B 2}
    The support of p is a compact set $\mathcal{X}$ such that, $p_{min} = \inf_{x \in \mathcal{X}} p(x) > 0$, and $p_{max} = \sup_{x \in \mathcal{X}} p(x)$.
\end{assumptionB}

\begin{assumptionB} \label{B 3}
    The response variable satisfies $Y_{min}< Y < Y_{max}$ almost surely, where $Y_{min}$ and $Y_{max}$ are real numbers.
\end{assumptionB}  

\begin{assumptionB} \label{B 4}
    The regression function $m$ is continuous, and $ \gamma < m(x) < M$ for all $x \in \mathbb{R}$, where $\gamma$ and $M$ are real numbers.
\end{assumptionB}  

\begin{assumptionB} \label{B 5}
    The derivative of the regression function $m_{1}$ is tame and continuous on a compact domain, and $\gamma_{1}< m_{1}(x) < M_{1}$ for all $x \in \mathbb{R}$, where $\gamma_{1}, M_{1} \in \mathbb{R}$.
\end{assumptionB}  

Now, we state the main result of the consistency of $\widehat{PH}^{\epsilon}_{k}(m_{1})$ in the following theorem.
\begin{t1} \label{Th:1}
 If $h \equiv h(n) \xrightarrow{} 0$ as $n \xrightarrow{} \infty$ such that $n h^6 \xrightarrow{} \infty$, then under the Assumptions A.~\ref{A 1} -- A.~\ref{A 3} and B.~\ref{B 1} -- B.~\ref{B 5}, we have   
\begin{equation} \label{Main Theorem}
    \mathbb{P}\left( \delta_{B} \left(\widehat{PH}_{k}^{\epsilon}\left(m_{1}\right), PH_{k}\left(m_{1}\right)\right) \leq 5\epsilon\right) \geq 1 - 3M_{\epsilon}n e^{-{C_{\epsilon/2}n h^6}},
\end{equation}
where k is any arbitrary degree of homology, $\delta_{B}$ is the bottleneck distance (see, Definition~\ref{Bottleneck distance}), $M_{\epsilon}$ is a positive integer defined as $M_{\epsilon} = \max \{i : t_i \geq \gamma_1\}$, $t_i$ is defined in Equation~\eqref{ Theorem constants},  and $C_{\epsilon}$ is defined as 
\begin{equation} \label{C_epsilon}
    C_{\epsilon} =  \displaystyle{\frac{8\delta^2 p_{min}^2 ( \gamma_1 + \epsilon)^2}{ (X_{max} - X_{min})^2 \left[\tau_2 (Y_{max} - Y_{min}\delta ) + \tau_1 Y_{min} \delta (1 - \delta ) \right]^2}} \text{, provided } \gamma_1 \neq -\epsilon.
\end{equation}
Here, all the involved constants are precisely defined in Assumptions A.~\ref{A 1} -- A.~\ref{A 3} and B.~\ref{B 1} --B.~\ref{B 5}   
\end{t1} 

In the following, we present a short sketch of the proof of Theorem \ref{Th:1}, while the full proof can be found in Section \ref{Proofs}. Let $\mathcal{D}$ denote the continuous filtration and $\mathcal{D}^{\epsilon}$ denote the discrete version of $\mathcal{D}$ with step size 2$\epsilon$. In the present context, we estimate the persistent homology corresponding to the filtration $\mathcal{D}^{\epsilon}$. Therefore, the maximum difference between the persistent homology of $\mathcal{D}$ and $\mathcal{D}^{\epsilon}$ is 2$\epsilon$. Hence, to prove Theorem \ref{Th:1}, it is enough to show that the maximum difference between estimated persistent homology and true persistent homology of $\mathcal{D}^{\epsilon}$ is 3$\epsilon$. As pointed out in \cite{Tdaconsistency}, this can be done using the language of $\epsilon$-interleaving introduced in \cite{Proximityofpd}. Thus, using the weak stability theorem from \cite{Proximityofpd}, and Lemma \ref{Lemma2}, the proof of Theorem \ref{Th:1} is established.
 
 \begin{r1} \label{Th1:remark}
Note that the constant $C_{\epsilon}$ must be positive in order to have a valid lower bound of the probability in Theorem \ref{Th:1}. Since if $C_{\epsilon} = 0$, then the sequence $1 - 3M_{\epsilon}n e^{-{C_{\epsilon/2}n h^6}} = 1 - 3M_{\epsilon}n < 0$. Furthermore, as $n \xrightarrow{} \infty$, the sequence goes to -$\infty$, implying that the proposed procedure is not consistent. For this reason, we need to assume that $\delta > 0$ as precisely defined in Assumption A.~\ref{A 2}.     
 \end{r1}

 \begin{r1} 
     As pointed out by one of the reviewers, it is noteworthy that the much slower contraction rate of the kernel bandwidth is a consequence of Assumption A.~\ref{A 3}. Precisely, Assumption A.~\ref{A 3} is used in Lemma \ref{Lemma1} to find bounds of the auxiliary random variables $Z_i$ (see Equation~\eqref{Zi}). Consequently, an extra factor of $h^2$ is introduced in the bounds on the auxiliary random variables $Z_i$. This leads to the contraction rate of the kernel bandwidth as $h \lesssim n^{-1/6}$ in Theorem \ref{Th:1}. Moreover, Assumption A.~\ref{A 3} explicitly ties the slow contraction rate of the kernel bandwidth to the fact that first derivatives are considered. 
 \end{r1}

\section{ Statistical significance of local structures} \label{Quantitative evaluation}
In the following sections, we apply Theorem \ref{Th:1} in identifying the local structures such as monotonicity, convexity, and modality in smooth regression curves. One major challenge in inferring local structures of regression functions from a finite data set is assessing their statistical significance. A quantitative statistical assessment is required to determine whether the observed local structures are signals or artifacts of sampling noise. In the present context, local structures are derived from the 0-dimensional estimated persistence barcodes computed using the procedure given in Section \ref{Problem Formulation and Main Results} whose validity is guaranteed by Theorem \ref{Th:1}. Therefore, we translate the statistical significance of the observed local structures in terms of the statistical significance of the 0-dimensional topological features (connected components) in the estimated barcode.     

There are various procedures to assess the statistical significance of topological features in estimated barcodes obtained from continuous filtrations by existing procedures (see, e.g., \cite{Fasy2014}). However, in the present context, we compute barcodes from a discrete filtration using the method proposed by \cite{Tdaconsistency} that does not have similar results as in the case of continuous filtration. The main reason is that most existing procedures for examining the statistical significance require stability results for persistence barcodes. In contrast, we do not have the stability results for persistence barcodes in the present context.   

We propose a measure of statistical significance of the observed 0-dimensional topological features in the estimated barcodes calculated from the procedure described in Section \ref{Problem Formulation and Main Results}. We provide a significance measure for the topological features in the observed barcode using Theorem \ref{Th:1}. For $\alpha \in (0, 1)$, a $(1-\alpha)$-confidence set using Theorem \ref{Th:1} can be constructed as: 
\begin{equation} \label{Condifence_set}
    \mathcal{C}_{n} := \left\{ \mathcal{B}: \delta_{B}\left(\mathcal{B}, \hat{\mathcal{B} }\right) \leq 5 z_{n, \alpha} \right\}, 
\end{equation}
where $\mathcal{B}$ is an arbitrary barcode, $\hat{\mathcal{B}}$ is a barcode estimated using the procedure given in Section~\ref{Problem Formulation and Main Results}, and $z_{n, \alpha}$ solves the following: 
\begin{equation} \label{Znalpha}
    3M_{z_{n, \alpha}}n \exp \left({-nh^6 C_{z_{n, \alpha}/2}}\right) = \alpha.
\end{equation}
Here $M_{z_{n, \alpha}}$ and $C_{z_{n, \alpha}}$ are defined by plugging $z_{n, \alpha}$ in place of $\epsilon$ in $M_{\epsilon}$ and $C_{\epsilon}$ (see Equation~\eqref{C_epsilon}), respectively. 

Note that $\mathcal{C}_{n}$ contains the true barcode $\mathcal{B}_0$ with probability $1 - \alpha$. That is
\begin{align*}
    \mathbb{P} \left( \mathcal{B}_0 \in \mathcal{C}_{n} \right) &= \mathbb{P} (\delta_{B}(\mathcal{B}_0, \hat{\mathcal{B} }) \leq 5 z_{n, \alpha})\\
    & = \mathbb{P} \left(\delta_{B} \left( PH_k(m_1), \widehat{PH}^{\epsilon}_k (m_1) \right) \leq 5z_{n, \alpha}\right)\\
    & \geq 1 - \alpha, \tag{By Theorem \ref{Th:1} and Equation \eqref{Znalpha}}
\end{align*}

Given a $(1 - \alpha)$-confidence set $\mathcal{C}_{n}$ for persistence barcodes, we can use the significance measure discussed in \cite{Fasy2014} to declare a bar to be significant. In particular, a bar $[a, b] \in \hat{\mathcal{B}}$ is significant at $\alpha$ level of significance, if  
\begin{equation} \label{significance measure}
    |a - b| > 5 \sqrt{2} z_{n, \alpha}.
\end{equation}

\subsection{Bandwidth Selection}\label{BS} 
We provide an appropriate choice of the kernel bandwidth $h\equiv h(n)$ to infer structures for a finite sample size $n$. Note that, for any arbitrary bandwidth $h$, the sample size $n$ must be large enough to ensure convergence of the estimation procedure from Theorem \ref{Th:1}. For example, for $h = 0.1$, we require $n > 10^6$ to ensure convergence from Theorem \ref{Th:1} provided that $C_{\epsilon/2}$ is not a sufficiently large quantity. This suggests choosing $h$ appropriately, depending on the sample size $n$ and the unknown population quantities involved in Theorem \ref{Th:1}. 

Given the unknown population quantities $M_\epsilon$ and $C_\epsilon$, we choose $h = h_l$ to infer structures at a $\alpha$ level of significance such that $h_l$ solves
\begin{equation} \label{Bandwidth_Selection}
    3M_{\epsilon}n \exp \left({-nh^6_l C_{\epsilon/2}}\right) = \alpha.
\end{equation}

In view of Equation~\eqref{Znalpha}, note that for all $h \geq h_l$, $z_{n, \alpha} = \epsilon$. In other words, by Theorem~\ref{Th:1}, $5\epsilon$ is $(1 - \alpha)$-th quantile of the sampling distribution of $\delta_{B}$ $ ( \widehat{PH}^{\epsilon}_k (m_1), PH_k(m_1))$, if the chosen bandwidth $h \geq h_l$. On the other hand, for $h \in (0, h_l)$, $z_{n, \alpha} > \epsilon$ such that $5z_{n, \alpha}$ is $(1 - \alpha)$-th quantile of the sampling distribution of $\delta_{B} ( \widehat{PH}^{\epsilon}_k (m_1), PH_k(m_1))$. This implies that the confidence set $\mathcal{C}_n$ is wider for $h \in (0, h_l)$, and sharper for $h \geq h_l$, for a finite sample size $n$. 

Now, we discuss the estimation of the unknown population quantities $M_{\epsilon}$ and $C_{\epsilon}$, from the observed data $(Y_1, X_1), \ldots, (Y_n, X_n)$. In this context, we discuss only the estimation of unknown quantities $\delta$, $p_{min}$, and $\gamma_{1}$ involved in $C_{\epsilon}$ (see Equation \eqref{C_epsilon}) and $M_{\epsilon}$, since the estimation of other unknown quantities is straightforward.   

It is straightforward to estimate $\delta $ from the chosen kernel as $\min K(X_i)$. However, $\min K(X_i)$ could be a very small quantity close to 0 as $K$ is bounded, which is not desirable as $\delta > 0$. In such a situation, we simply take a reasonably larger value than 0, say $\delta = 0.1$. Therefore, in what follows, we take $\delta = \max (0.1, \min K(X_i) ) $ to perform data analysis. One can use kernel density estimates of the regressors to estimate $p_{min}$. However, due to kernel artifacts at boundaries, the estimate of $p_{min}$ is close to 0. Note that $p_{min} > 0$, and the smaller value will make $C_{\epsilon}$ close to 0, which obstructs the validity of Theorem \ref{Th:1}. Therefore, we assume that $p_{min} \geq 0.1$ and use this lower bound, which will be sufficient for our purpose.

Next, to estimate $\gamma_{1}$ and $M_{\epsilon}$, we propose Algorithm~\ref{Algorithm}. Consider the definition of derivative for $x_1 < x_2 \in \mathbb{R}$:
\begin{equation} \label{Derivative_def}
    m_1(x_1) = \lim_{x_2\downarrow x_1} \frac{m(x_2) - m(x_1)}{ x_2 - x_1} \approx (m(x^*) - m(x_1)) \eta^{-1}, 
\end{equation}
where $\eta > 0 $ such that $x_2 - x_1 \approx \eta$, and $x^* \in (x_1, x_1 + \eta)$.

Then, given the noisy observations $Y_1, \ldots, Y_n$ of $m$ at the points $X_1, \ldots, X_n$, respectively. The procedure of estimating derivative values is articulated in Algorithm~\ref{Algorithm}.
\begin{algorithm} 

\caption{Derivative Estimation Procedure}
\label{Algorithm}

\begin{algorithmic}[1] 
\State \textbf{Input:} Data $\{(Y_1, X_1), \ldots, (Y_n, X_n)\}$.
\State \textbf{Step 1: Sorting.} 
Sort the regressors in increasing order $X_{(1)}, \dots, X_{(n)}$. 
\State \textbf{Step 2: Threshold Selection.}
Define the threshold:

$\eta = Median \{X_{(2)} - X_{(1)}, \ldots, X_{(n)} - X_{(n - 1)}  \}$ 

\State\textbf{Step 3: Filtration.} 
Filter $X_i$'s and $Y_i$'s as follows:

$X_\eta = \{ X_{(i)} : X_{(i + 1)} - X_{(i)} < \eta, i = 2, \ldots, n\}$, and 
$ Y_\eta = \{ Y_{(i)}: X_{(i)} \in X_\eta\}$.

\State \textbf{Step 4: First-Order Differences.} 
For points in $X_\eta$ and $Y_{\eta}$, compute:

$D_j = \frac{Y_{(j+1)} - Y_{(j)}}{X_{(j+1)} - X_{(j)}}, \quad j = 1, \dots, N$, $N = |X_{\eta}| = |Y_{\eta}|$. 

\State \textbf{Step 5: Cleaning.} Mitigate noise in $D_j, j =1, \ldots, N $ by defining:

$\hat{Y}^{(1)} = \{ D_j: D_j \in [Q(0.25) - 1.5IQR, Q(0.75) + 1.5IQR]\}$, 

where $Q(p)$ denotes $p$-th empirical quantile of $\{D_j: j = 1, \ldots, N\}$, and $IQR = Q(0.75) - Q(0.25)$.  

\State \textbf{Output:} Derivative estimates $\hat{Y}^{(1)}$. 
\end{algorithmic}
\end{algorithm}

Now, using the Algorithm~\ref{Algorithm}, the estimates of $\gamma_1$ and $M_{\epsilon}$ are defined as:
\begin{equation}
   \hat{\gamma}_{1} := \min_i \hat{Y}^{(1)}_i  \text{ and } \hat{M}_{\epsilon} := \left \lceil  \hat{N}_{\epsilon} - \frac{\hat{\gamma}_{1} }{2\epsilon} \right \rceil  \text{, where } \hat{N}_{\epsilon} := \max_i \lceil \hat{Y}^{(1)}_i / 2 \epsilon \rceil, \hat{Y}^{(1)}_i \in \hat{Y}^{(1)}.
\end{equation}
Note that using derivative estimates as $\frac{Y_i - Y_{i - 1}}{X_i - X_{i-1}}$ is too noisy in general (see, e.g., \cite{DerLPR}) that may affect the results adversarially. 

\begin{r1} \label{remark:Algo}
    The main idea of Algorithm~\ref{Algorithm} is motivated by Equation \eqref{Derivative_def}, which suggests approximating the true derivative only for those points in the domain for which there exists $x^*$ in a small neighborhood. In essence, this suggests avoiding approximating $m_1$ in sparse regions of the domain for a good approximation by Equation \eqref{Derivative_def}. We avoid these sparse regions by defining a suitable threshold $\eta$ in Step 2, where we have taken $\eta = Median\{X_{(2)} - X_{(1)}, \ldots, X_{(n)} - X_{(n - 1)}\}$. However, $\eta$ can be suitably chosen based on the data such that $\eta > 0$ and $X_{\eta} \neq \varnothing$ in Step 3. In general, $\eta = Median\{X_{(2)} - X_{(1)}, \ldots, X_{(n)} - X_{(n - 1)}\}$ may not be a good choice, as it may be 0 for certain data sets. For instance, in Section~\ref{Data set 1} for cars data,  $Median\{X_{(2)} - X_{(1)}, \ldots, X_{(n)} - X_{(n - 1)}\} = 0$. Therefore, we define $\eta = mean \{X_{(2)} - X_{(1)}, \ldots, X_{(n)} - X_{(n - 1)}\}$ for this data set which is > 0 and works perfectly.

    \noindent Note that $X_\eta$ in Step 3 provides suitable values in the domain from the observed data for which $X_{(i + 1)}$ plays the same role for $X_{(i)}$ as $x^*$ plays for $x_1$ in Equation \eqref{Derivative_def}, for all $i = 2, \ldots, n$. Thus, for a suitably chosen $\eta$, the ratio of first-order differences of $Y \in Y_{\eta}$ (see Step 3) and $X \in X_\eta$ (see Step 3) provides good estimated values of $m_1$ by Equation \eqref{Derivative_def}, at the points in $X \in X_\eta$. However, due to randomness in $X_{i}$'s, these estimated values can still be noisy. Therefore, we mitigate noise in the estimated derivative values by cleaning them in Step 5.  
\end{r1}

We implement Algorithm~\ref{Algorithm} in simulation and real data analysis to estimate $\gamma_1$ and $M_{\epsilon}$ in the following sections. It is evident from simulation and real data analysis results in Section~\ref{Simulation studies} and Section~\ref{Real data study} that Algorithm~\ref{Algorithm} works well to infer true structures underlying the data based on Theorem \ref{Th:1} and Equation \eqref{Znalpha} or Equation \eqref{Bandwidth_Selection}. Moreover, we provide an explicit implementation of Algorithm~\ref{Algorithm} in Appendix~\ref{Alogorithm 1} for the functions used in Section~\ref{Simulation studies} to examine whether it provides reasonably good estimates of the true derivative.  

\begin{r1}
    It is an appropriate place to mention that one of the reviewers suggested the following modifications to Algorithm~\ref{Algorithm} to circumvent the issue that leads to resorting to any other threshold than the median in Step 2 of Algorithm~\ref{Algorithm}. The reviewer suggested a preprocessing step at the start of Algorithm~\ref{Algorithm}, which eliminates all the data points $(X_i, Y_i)$ for which there is a $j \neq i$ such that $X_j = X_i$ or aggregates all such points into one point $(X, \bar{Y})$, where $\bar{Y}:= mean \{Y_i: i \in \mathcal{I}\}$, $\mathcal{I}$ is a finite index set such that for any $i, j \in \mathcal{I}$, we have $X_i= X_j$, and $X := X_i$ for all $i \in \mathcal{I}$. This preprocessing step would ensure that the difference quotient $D_j$ in Step 4 is not ill-defined. Moreover, it also ensures that for any two consecutive points $X_{(i)}$ and $X_{(i+1)}$ such that either $X_{(i)}$ or $X_{(i + 1)}$ is non-unique, the corresponding value of the regression functions is well defined. As such, the proposed preprocessing step would rule out the possibility of threshold $\eta = 0$ in Step 2 of Algorithm~\ref{Algorithm}. This makes the median viable again and avoids resorting to the mean as discussed in the remark after Algorithm~\ref{Algorithm}. The aforementioned preprocessing step can be useful in practice for many datasets, such as Dataset 1 in Section \ref{Data set 1}. However, it is not used in this article; instead, in view of the fact explained in the first paragraph in Section \ref{Data set 1}, a mean-based threshold is used in Step 2 of Algorithm~\ref{Algorithm} for Dataset 1 studied in Section \ref{Data set 1}.    
\end{r1}

\section{Simulation studies} \label{Simulation studies}
 This section presents some applications of Theorem \ref{Th:1} to infer the monotonicity, convexity, and modality of smooth regression curves underlying the data. We simulate noisy observations from smooth regression curves and measure the statistical significance of the observed local structures using the significance measure proposed in Section \ref{Quantitative evaluation}. 

 In addition, we compare the results with a graphical device called SiZer map, an existing procedure for exploring structures such as peaks and valleys in smooth regression curves proposed by \cite{SiZer}. The SiZer map visualizes significant features such as bumps simultaneously over both location and scale (bandwidth) by using a color map. The idea is to highlight the regions in the location-scale space where the curve significantly increases and decreases. The idea is motivated by the fact that the significant bumps will be at zero crossings of the derivative between regions of significant increase and decrease. Hence, the name ``SiZer" is given for ``SIgnificant ZERo crossings of derivatives". The color scheme on the SiZer map is as follows: \textbf{Blue} is used to indicate regions where the curve is significantly increasing; \textbf{Red} is used to indicate regions where the curve is significantly decreasing; \textbf{Purple} is used to indicate regions where the curve cannot be concluded to be either significantly increasing or decreasing; \textbf{Gray} is used to indicate regions where there are not enough data points to make statements about the significance of the features. The dotted curves in SiZer maps show effective window widths for each bandwidth as intervals representing $\pm2h$. 

In this section, we simulate noisy observations with the noise distribution as the truncated normal distribution supported on [-1, 1] with mean 0 and standard deviation 0.01. We estimate the regression function using the truncated standard normal and truncated standard Cauchy kernel supported on [-1, 1]. In particular, both these kernels take the form $K(x):= g(x)/g(0)$ to satisfy Assumptions A.~\ref{A 1}, A.~\ref{A 2}, and A.~\ref{A 3}, where $g$ is taken here to be the probability density function of either the truncated standard normal or the truncated standard Cauchy distribution supported on [-1, 1]. That is, $g(x) := f(x) /\{F(1) - F(-1)\}, x \in [-1, 1]$, where $f$ and $F$ are the probability density and the distribution function of either the standard normal or the standard Cauchy distribution, respectively. Moreover, it is an appropriate place to mention that we compare the results with the SiZer maps only for the truncated standard normal supported on [-1, 1], since the conventional SiZer is applicable only for the Gaussian kernel.   

\subsection{ Monotonicity} \label{ Monotonicity} 

We characterize the monotonicity of a smooth and bounded regression function using the 0-dimensional features of the barcode of the first derivative. The main idea is that the birth and death times corresponding to the 0-dimensional feature with the largest persistence in a barcode provide the maximum and minimum value of the function, respectively (see Figure~\ref{fig: persistent homology}). Therefore, we can extract the 0-dimensional feature with the largest persistence from the barcode of the first derivative to characterize the monotonicity of the function. In particular, if the death time of the feature is non-negative, that is, the minimum value of the first derivative is non-negative, then the function is non-decreasing. Similarly, if the birth time of the feature is non-positive, that is, the maximum value of the first derivative is non-positive, then the function is non-increasing.

Consider the regression function $\displaystyle{m(x) = e^{x}, x \in} $ [-1, 1] to investigate its monotonicity from the kernel estimates of its derivative $m_{1}$. We obtain barcodes for the estimated derivative, $\hat{m}_{1}$, using the procedure mentioned in Section \ref{Problem Formulation and Main Results}. The barcodes for the truncated standard Gaussian kernel are $\boldsymbol{\{ [0.13, 0.14]\}}$, $ \boldsymbol{\{ [0.15, 0.17]\}}$ and $ \boldsymbol{\{ [0.17, 0.20]\}}$ for the bandwidths $h = 1.5~(n = 200)$, $h = 1.37~(n = 400)$ \text{ and } $h = 1.29~(n = 600)$, respectively. These barcodes are obtained for the filtration parameter $\epsilon = 0.001$ and the bandwidths are chosen according to Equation \eqref{Bandwidth_Selection} to assess the statistical significance of observed barcodes at $\alpha = 0.05$ level of significance. Hence, by using Equation \eqref{significance measure} for $z_{n, \alpha} = 0.001$, we conclude that each of the observed bars is significant for the Gaussian kernel. Since the length of each bar (persistence) is greater than $5 \sqrt{2} z_{n, \alpha} = \boldsymbol{0.007}$. This implies that the observed data suggest, with the 95\% level of confidence, that the underlying function is strictly increasing, as the death times are positive for each of the bars.      

\begin{figure} [htbp!]
    \centering
    \includegraphics[width=\linewidth]{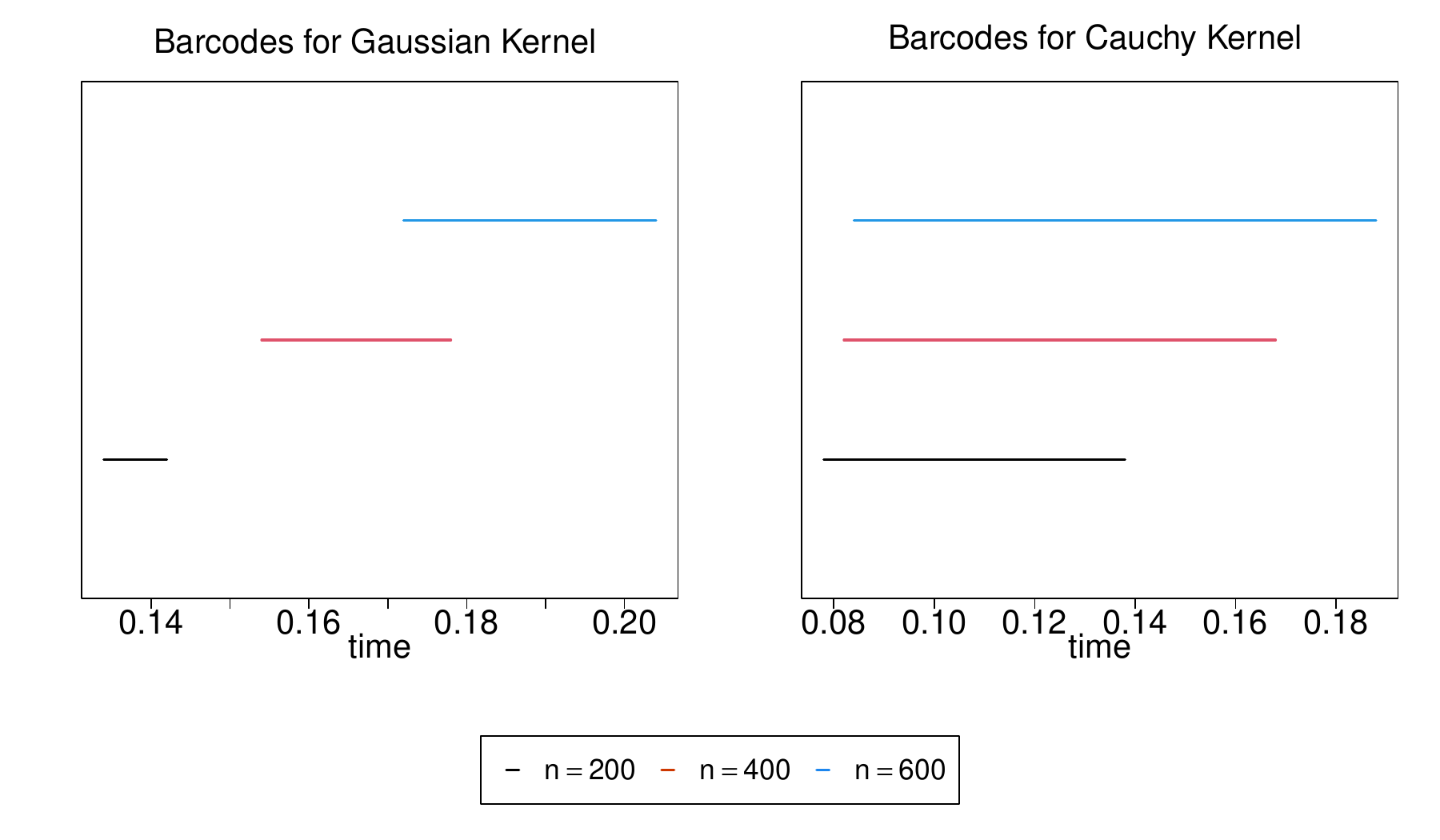}
    \caption{Estimated barcodes of $m_1$ obtained from noisy data simulated from the function $m(x) = e^x, x\in[-1, 1]$ (see Figures \ref{plot_combined} and \ref{Cauchy_combine}).}
    \label{n_50}
\end{figure}

The same conclusion also holds for the truncated standard Cauchy kernel. Since, the barcodes for the Cauchy kernel are $\boldsymbol{\{ [0.07, 0.13]\}}$, $ \boldsymbol{\{ [0.08, 0.16]\}}$ and $ \boldsymbol{\{ [0.08, 0.18]\}}$ for the bandwidths $h = 2~(n = 200), h = 1.8~(n = 400) \text{ and } h = 1.7~(n = 600)$, respectively. Therefore, by Equation \eqref{significance measure}, all the observed bars are significant. Note that, here also, we have taken $\epsilon = 0.001$ and the bandwidths are obtained using Equation \eqref{Bandwidth_Selection} for $\alpha = 0.05$. The estimated barcodes for both kernels are shown in Figure \ref{n_50}. 

\begin{figure} [htbp!]
    \centering
    \includegraphics[width=\linewidth]{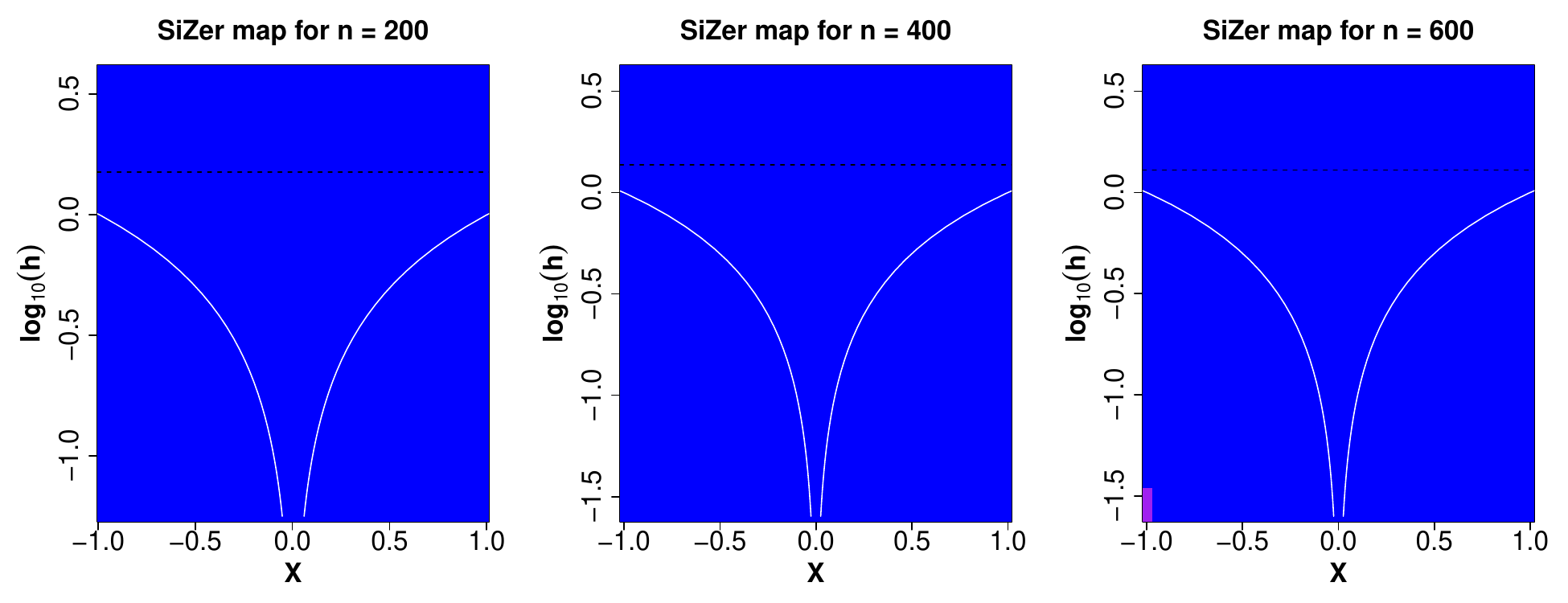}
    \caption{SiZer analysis of data (see Figures \ref{plot_combined} and \ref{Cauchy_combine}), to investigate the monotonicity of the underlying true function $m(x) = e^x, x\in[-1, 1]$.}
    \label{SiZer_monotonicity}
\end{figure}

Now, we compare the inferred structure from the barcode of $\hat{m}_1$ with the SiZer map of the data. It is evident from the SiZer map shown in Figure \ref{SiZer_monotonicity} that the observed structure is strictly increasing for each of the bandwidths used for the Gaussian kernel. Since the color of the scale-space region for each of the bandwidths (highlighted by the black dotted lines) is blue. 

\begin{figure}[htbp!]
    \centering
    \includegraphics[width=\linewidth]{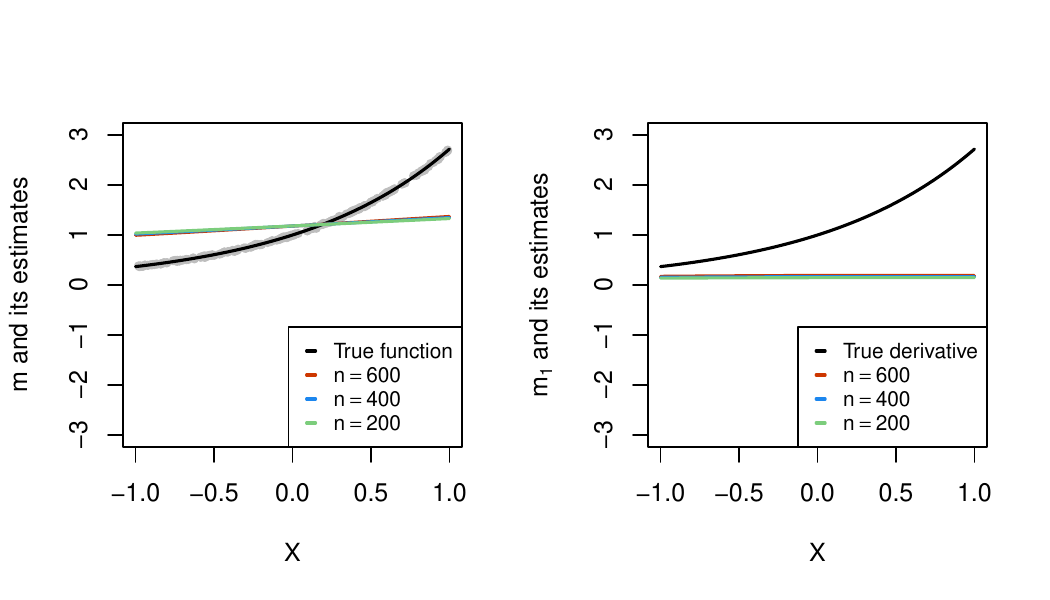}
    \caption{The true regression $m$ and its derivative $m_1$, $m_1(x)= m(x) = e^x, x\in [-1, 1]$ and their kernel estimates of for $h = 1.5~(n = 200), 1.37~(n = 400)\text{ and } 1.29~(n = 600) $, using the truncated Gaussian kernel on [-1, 1].}
    \label{plot_combined}
\end{figure}

\begin{figure} [htbp!]
    \centering
    \includegraphics[width=\linewidth]{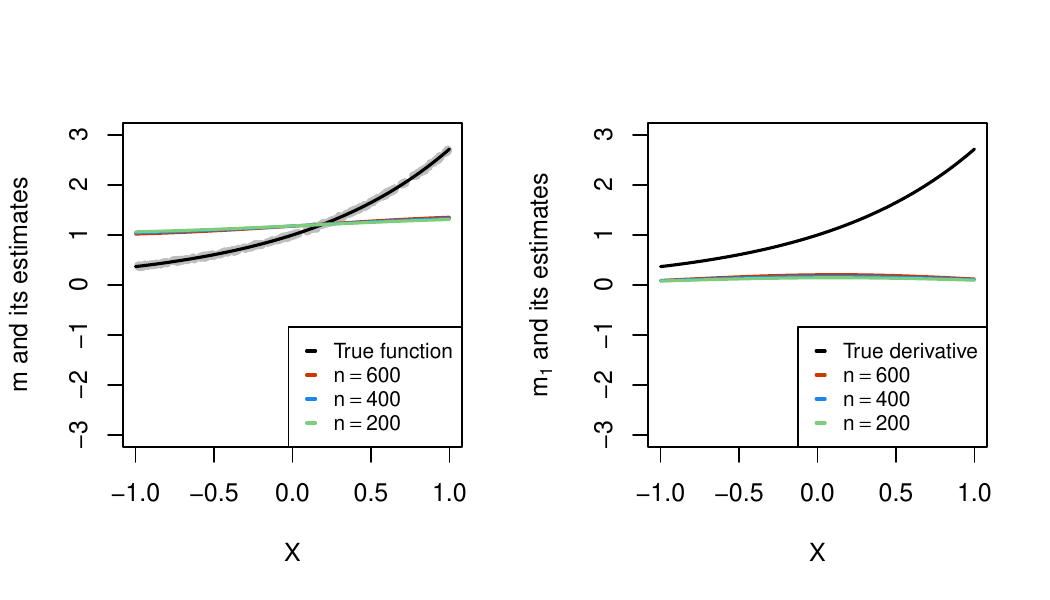}
    \caption{The true regression $m$ and its derivative $m_1$, $m_1(x)= m(x) = e^x, x\in [-1, 1]$ and their kernel estimates of for $h = 2~(n = 200), 1.8~(n = 400)\text{ and } 1.7~(n = 600) $, using the truncated Cauchy kernel on [-1, 1].}
    \label{Cauchy_combine}
\end{figure}

Note that we do not use point-wise or uniform consistency of the estimators of $m$ or $m_1$ here. It is evident from Figure \ref{plot_combined} and Figure \ref{Cauchy_combine} that neither $m$ nor $m_1$ is well estimated for each of the chosen bandwidths. Furthermore, the monotonicity (non-decreasingness) of the underlying regression function $m$ is not evident from the visual inspection of the estimates of $m$. However, we have shown that the local structures in $m$ can be inferred using the barcodes of $m_1$ without using pointwise or uniform consistency of $m$ or $m_1$. This is because Theorem~\ref{Th:1} establishes the consistency of the estimated persistent homology of $m_1$ without stability results (see Theorem \ref{stability theorem}).
%%%%%%%%%%%%%%%%%%%%%%%%%%%%%%%%%%%%%%%%%%%%%%%%%%%%
\subsection{Convexity} \label{Convexity} 
We provide an illustrative example to infer the convexity of a smooth regression function from noisy samples using estimated barcodes of its first derivative. We characterize the convexity (equivalently, concavity) of a function $f$ on [-1, 1] by the fact that there exists a critical point of $f$, say, $c \in$ [-1, 1] such that $f$ is non-decreasing (non-increasing) in $[c, 1]$ and non-increasing (non-decreasing) in $[-1, c]$. However, as one of the reviewers pointed out, this characterization is not valid for some functions. For instance, consider the function $\displaystyle{m(x) = \frac{x^2}{1 + x^2 + x|x|}}, x \in $[-1, 1]. The reason is that the convexity of $m$ requires its first derivative $m_{1}$ to be non-decreasing in [-1, 1], while the characterization captures only the monotone structure of $m$ on $[c, 1]$ and $[-1, c]$ rather than the monotonicity of $m_{1}$ on [-1, 1]. This suggests modifying the characterization to account for the monotone structure of $m_{1}$ on [-1, 1] to infer the convexity of $m$ on [-1, 1].  

\noindent Note that to characterize the monotone structure of $m_{1}$ quantitatively, we need to consider the barcodes of the second derivative of $m$. However, qualitatively, the monotonicity of $m_{1}$ can be captured by comparing the barcodes of $m_{1}$ and -$m_{1}$, as also suggested by one of the reviewers. The rationale behind the argument is that for a monotone function $f$, the super-level set $D_{t} = \{x : f(x) \geq t\}$ of $f$ has only a single connected component, for every $t \in \mathbb{R}$. However, this is not a necessary condition since if $f$ is a concave shape function, then also $D_{t}$ will have only a single connected component for every $t \in \mathbb{R}$. Therefore, to rule out such situations, we have to make use of the property of monotone functions that if $f$ is monotone increasing (decreasing) on an interval $[a, b]$ (here we take $a$ = -1, $b$ =1), then -$f$ will be monotone decreasing (increasing) on $[a, b]$. This property allows us to characterize the monotonicity of a function $f$ from the barcodes of $f$ and -$f$. Note that in Section~\ref{ Monotonicity}, we characterize the monotonicity quantitatively using the barcode of the derivative. 

Therefore, we modify the earlier characterization and propose the following two-step procedure to characterize the convexity (equivalently, concavity) of a function $f$ on [-1, 1] using its first derivative $f^{'}$:

 \begin{enumerate} [label=(\roman*)]
    \item First, compare the barcode of $f^{'}$ and -$f^{'}$, qualitatively. In particular, compare the number of 0-dimensional features in the barcode of $f^{'}$ and -$f^{'}$. If there is more than one 0-dimensional feature in the barcode of $f^{'}$ or -$f^{'}$, conclude that $f$ is not convex (or concave). Otherwise, proceed with the following step.

    \item Characterize the convexity (concavity) of a function $f$ on [-1, 1] by the fact that there exists a critical point of $f$, say, $c\in$ [-1, 1] such that $f$ is non-decreasing (non-increasing) on $[c, 1]$ and non-increasing (non-decreasing) on $[-1, c]$.
\end{enumerate}

\begin{figure}[htbp!]
    \centering
    \includegraphics[width=\linewidth]{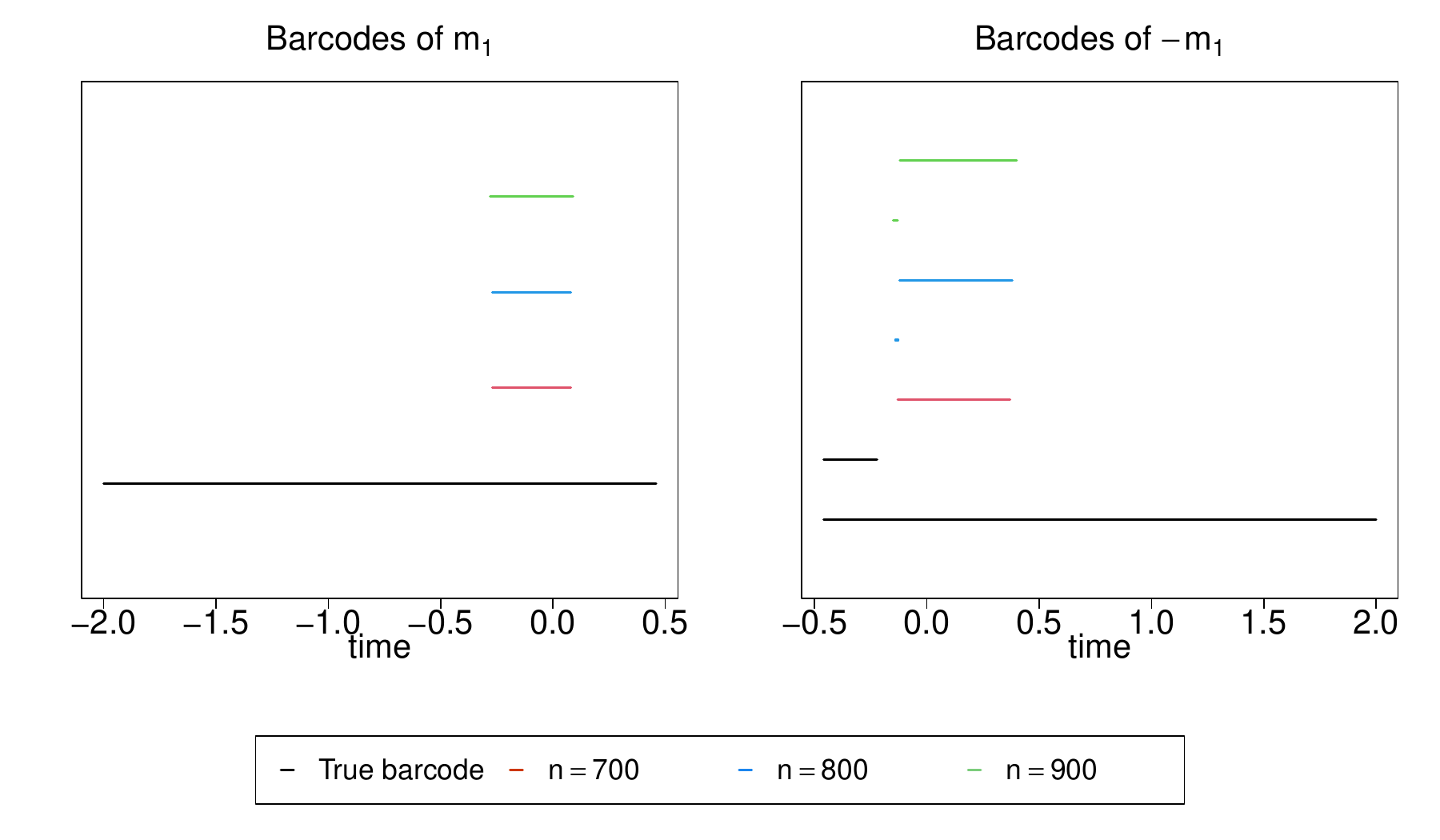}
    \caption{The barcode of $m_{1}$ and -$m_{1}$ along with the estimated barcodes obtained from noisy data simulated from the function $m(x) = x^2(1+ x^2 + x|x|)^{-1}, x\in [-1, 1]$ (see Figure \ref{convex function}), using the truncated Gaussian kernel on [-1, 1].} 
    \label{convexity_m}
\end{figure}

In the above characterization, step (i) ensures that $f^{'}$ is monotone on [-1, 1] or not. Since any structure impeding the monotone structure of the derivative would result in the birth of new connected components in the barcodes of either $f^{'}$ or -$f^{'}$, making a qualitative difference in the barcodes of $f^{'}$ and -$f^{'}$. This fact is evident from Figure \ref{convexity_m} where the true barcode of $m_{1}$ and -$m_{1}$ are qualitatively different. On the other hand, step (ii) is required to characterize the convexity (concavity) without using the second derivative of the function. We use the characterization of monotonicity proposed in Section \ref{ Monotonicity} to investigate the monotonicity of $f$ on $[-1, c]$ and $[c, 1]$ in step (ii).    

We infer the convexity of the function $m$ aforementioned in this subsection and compare the results with the respective SiZer maps. To infer the convexity of $m$ from the observed data, we need to consider the significant features in the observed barcode to apply step (i) of the proposed characterization. 
Note that, since we aim to infer the convexity of $m$ on the entire domain [-1, 1], the 0-dimensional feature with the largest persistence must be significant to proceed with step (ii), else we conclude that the observed structure is not convex (concave) significantly. The reason is that the features with shorter persistence represent structures in the subset of the domain [-1, 1], whereas here we need to infer the monotonicity of $m_{1}$ on [-1, 1]. 

We simulate noisy observations from $m$ and estimate $m_1$ for the bandwidths chosen according to Equation \eqref{Bandwidth_Selection} for $\alpha = 0.05$, and obtain barcodes for $\epsilon = 0.001$. The estimated barcodes of $m_1$ using the truncated Gaussian kernel are $\boldsymbol{\{ [-0.27, 0.08]\}}$, $\boldsymbol{\{ [-0.27, 0.08]\}}$ and $\boldsymbol{\{ [-0.28, 0.09]\}}$ for $h = 0.57~(n = 700)$, $h = 0.56~(n = 800)$ and $h = 0.55~(n = 900)$, respectively. Next, we estimate barcodes of -$m_1$ to investigate the convexity of $m$ using step (i) of the proposed characterization. The estimated barcodes of -$m_1$ for the same aforementioned parameters are $\boldsymbol{\{ [-0.13, 0.37]\}}$, $\boldsymbol{\{ [-0.14, 0.38], [-0.128, -0.122]\}}$ and $\boldsymbol{\{ [-0.15, 0.40], [-0.13, -0.12]\}}$ for $h = 0.46~(n = 700)$, $h = 0.45~(n = 800)$ and $h = 0.44~(n = 900)$, respectively. The true and estimated barcodes for $m$ and -$m_1$ are shown in Figure \ref{convexity_m}.

It is evident from Figure \ref{convexity_m} that as $n$ increases, the true features of $m$ and -$m_1$ are recovered. The statistical significance of the observed bars is evident from Equation \eqref{significance measure} for $z_{n, \alpha} = 0.001$, as the chosen bandwidths are the solution of Equation \eqref{Bandwidth_Selection} for $\alpha = 0.05$ and $\epsilon = 0.001$. Note that the second bar in the barcode of -$m_1$ is not significant, as its persistence is less than 0.007. However, as we increase $n$ from 800 to 900, the barcode of -$m_1$ consists of two significant bars. Hence, we conclude that the observed data suggest that the underlying function $m$ is not convex with the 95\% level of confidence.  

\begin{figure}[htbp!]
    \centering
    \includegraphics[width=\linewidth]{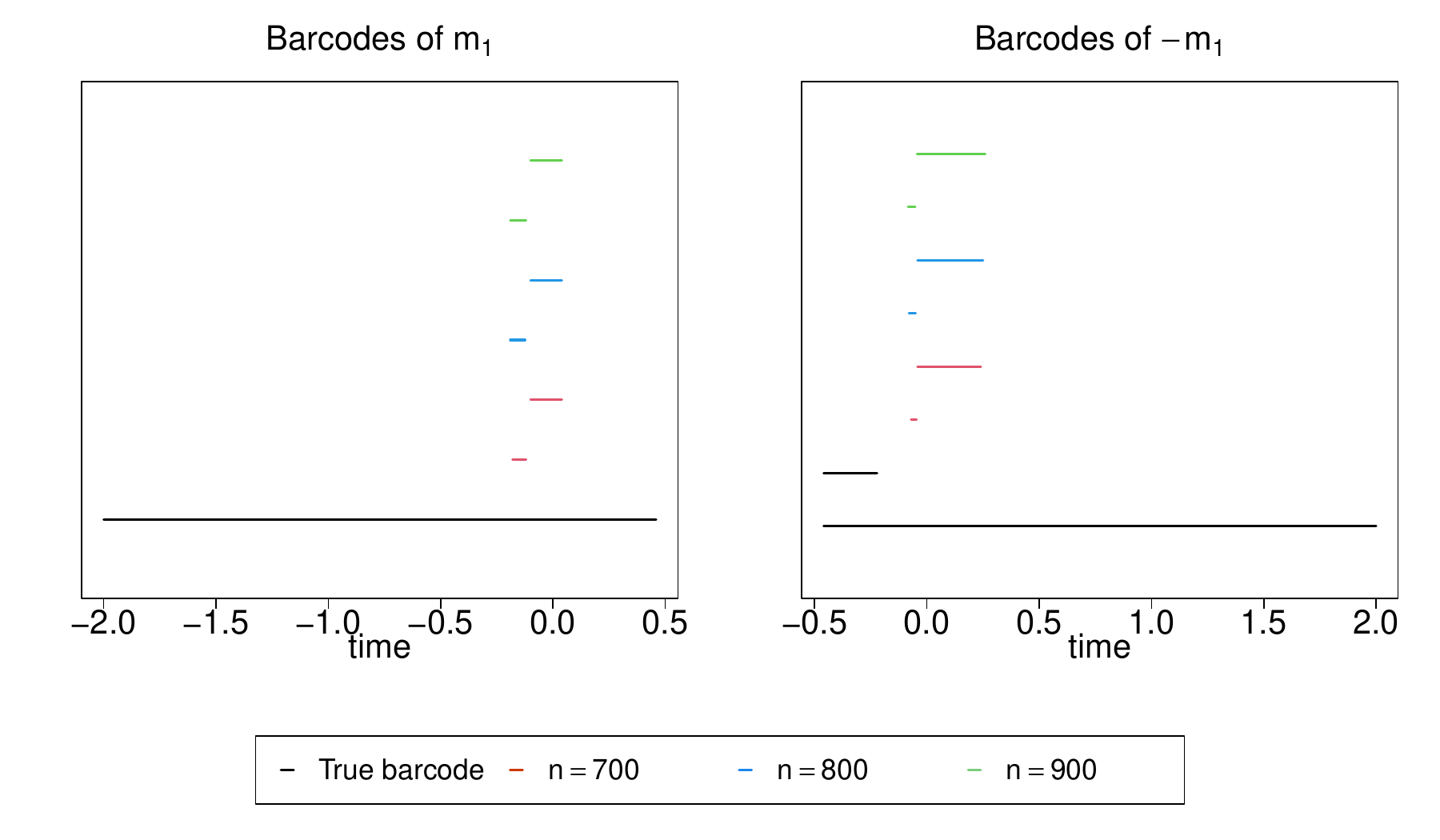}
    \caption{The barcode of $m_{1}$ and -$m_{1}$ along with the estimated barcodes obtained from noisy data simulated from the function $m(x) = x^2(1+ x^2 + x|x|)^{-1}, x\in [-1, 1]$ (see Figure \ref{convex function}), using the truncated Cauchy kernel on [-1, 1].}
    \label{Cauchy_convexity}
\end{figure}

Alternatively, we obtain barcodes of $m_1$ using the truncated Cauchy kernel to infer the convexity of $m$ from the observed data. The estimated barcodes of $m_1$ are $\boldsymbol{\{[-0.18, 0.04],[-0.12, -0.10] \}}$, $\boldsymbol{\{ [-0.19, 0.04], [-0.12, -0.10]\}}$ and$\boldsymbol{\{ [-0.19, 0.04], [-0.12, -0.10]\}}$ for $h = 0.77~(n = 700)$, $h = 0.75~(n = 800)$ and $h = 0.73~(n = 900)$ , respectively. Subsequently, to proceed with step (i) of the proposed characterization, we estimate barcodes of -$m_1$. The barcodes of -$m_1$ are $\boldsymbol{\{ [-0.07, 0.24], [-0.046, -0.042]\}}$, $\boldsymbol{\{ [-0.08, 0.25], [-0.05, -0.04]\}}$ and $\boldsymbol{\{ [-0.08, 0.26], [-0.05, -0.04]\}}$ for $h = 0.62~(n = 700)$, $h = 0.61~(n = 800)$ and $h = 0.6~(n = 900)$, respectively. These barcodes are shown in Figure \ref{Cauchy_convexity}. 

It is evident from the significance measure articulated in Equation \eqref{significance measure} that all the observed bars are significant except for the second bar in the observed barcode of -$m_1$ for $n = 700$. Note that here the significance measure is used for the bandwidths that solve Equation \eqref{Bandwidth_Selection} for $\alpha = 0.05$ and $\epsilon = 0.001$. Therefore, a bar with a length less than 0.007 is considered to be insignificant. Hence, by step (i), the observed data suggest that $m$ is not convex with the 95\% level of confidence, as there are two significant 0-dimensional features in the barcode of both $m_1$ and -$m_1$ as $n$ increases. 

\begin{figure}[htbp!] 
    \centering
    \includegraphics[width=\linewidth]{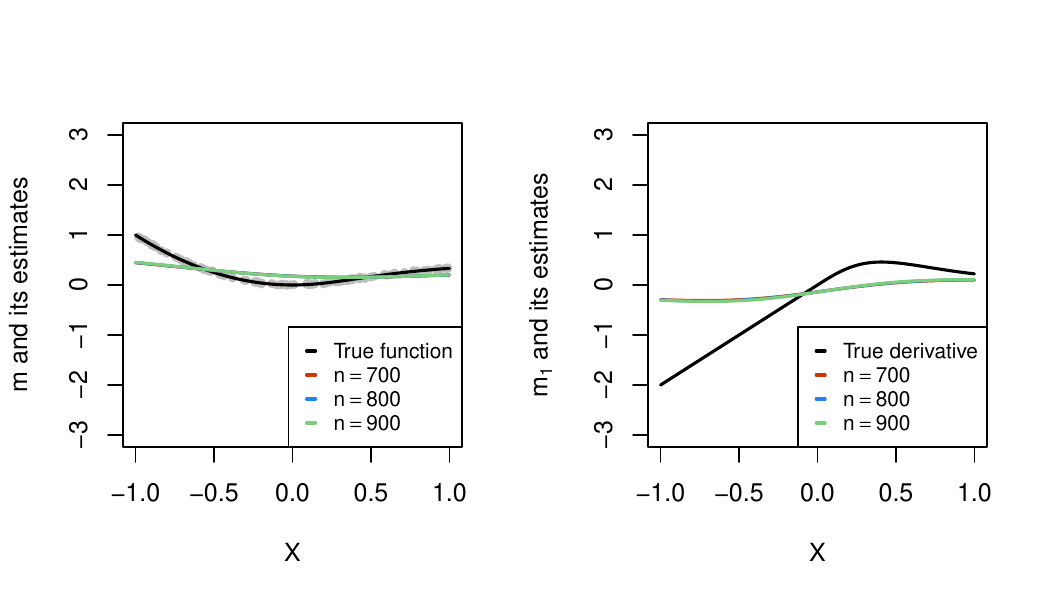}
    \caption{The true regression function $m(x) = x^2(1+ x^2 + x|x|)^{-1}, x\in [-1, 1]$, its derivative $m_1$ and their kernel estimates using the truncated Gaussian kernel on [-1, 1] for $h = 0.57~(n = 700)$, $h = 0.56~(n = 800)$ and $h = 0.55~(n = 900)$.}
    \label{convex function}
\end{figure}

Thus, Theorem \ref{Th:1} provides a valid inference on the observed structure using estimates of $m$ and $m_1$ that are pointwise not well estimated, as can be seen in Figure \ref{convex function}. Moreover, the visual inspection of the estimated curves in Figure \ref{convex function} does not reveal the true structure underlying the data.  

\begin{figure}[htbp!]
    \centering
    \includegraphics[width=\linewidth]{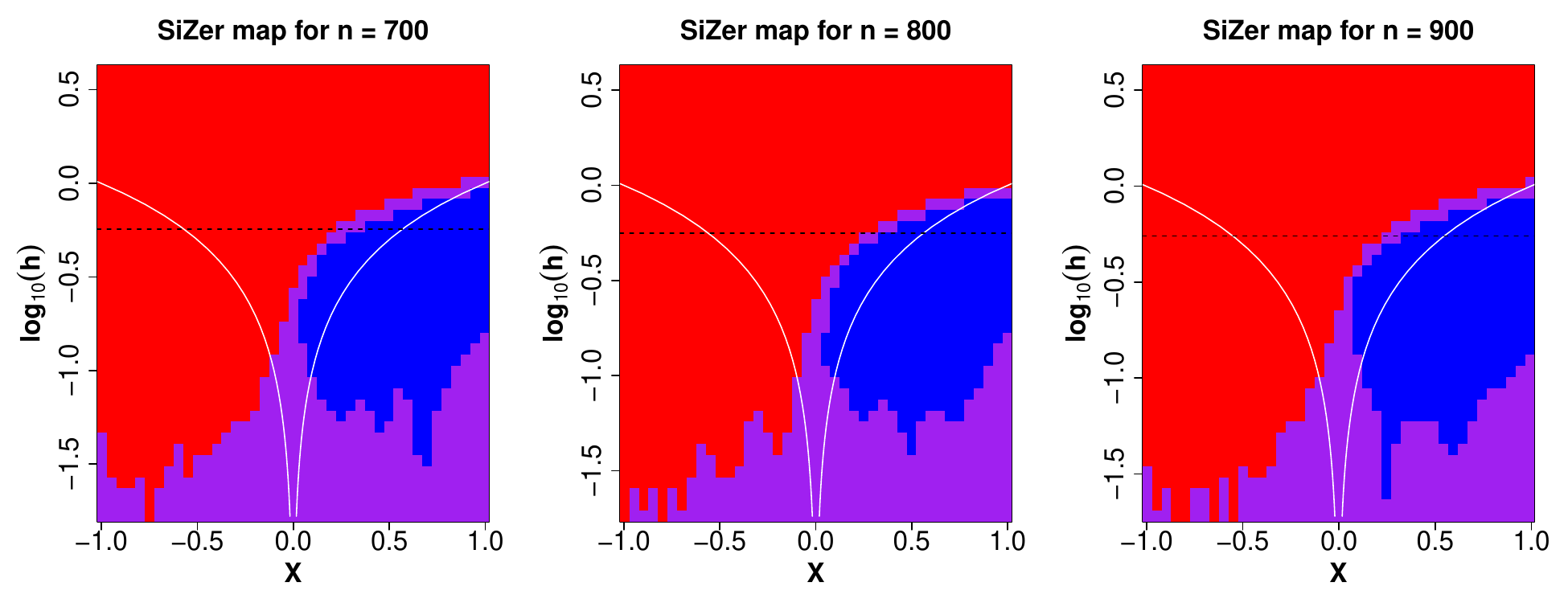}
    \caption{SiZer analysis of data shown in Figure \ref{convex function} to investigate the convexity of the underlying true function $m(x) = x^2(1+ x^2 + x|x|)^{-1}, x\in [-1, 1]$. }
    \label{SiZer_convexity}
\end{figure}

Now, we compare the results with the SiZer map of the data shown in Figure \ref{SiZer_convexity}. The SiZer maps in Figure \ref{SiZer_convexity} are inconclusive regarding the convexity of the observed structure for the bandwidths highlighted by the black dotted lines on the SiZer maps that are chosen according to Equation \eqref{Bandwidth_Selection} for the Gaussian kernel. This is because the SiZer maps merely detect a significantly decreasing slope followed by a significantly increasing slope of the estimated function, but SiZer does not have the suitable tools to conclude anything concerning convexity. However, the proposed methodology correctly infers the true structure of the underlying regression function $m$.      

%%%%%%%%%%%%%%%%%%%%%%%%%%%%%%%%%%%%%%%%%%%%%%%%%%%%%%%%%%%%%%%%%%
\subsection{Modality} \label{Modality}
We provide an illustrative example to infer the modality of a smooth regression function from noisy samples using estimated barcodes of its first derivative. We characterize the modes of a smooth regression function $f$ in [-1, 1] using the critical points of $f$. We say that a critical point $c \in $ [-1, 1] is a mode of $f$ if there exists a $\delta > 0$ such that $f$ is non-decreasing on $[c - \delta, c]$, and non-increasing on $[c, c + \delta]$. 

\begin{figure}[ht!]
    \centering
    \includegraphics[width=\linewidth]{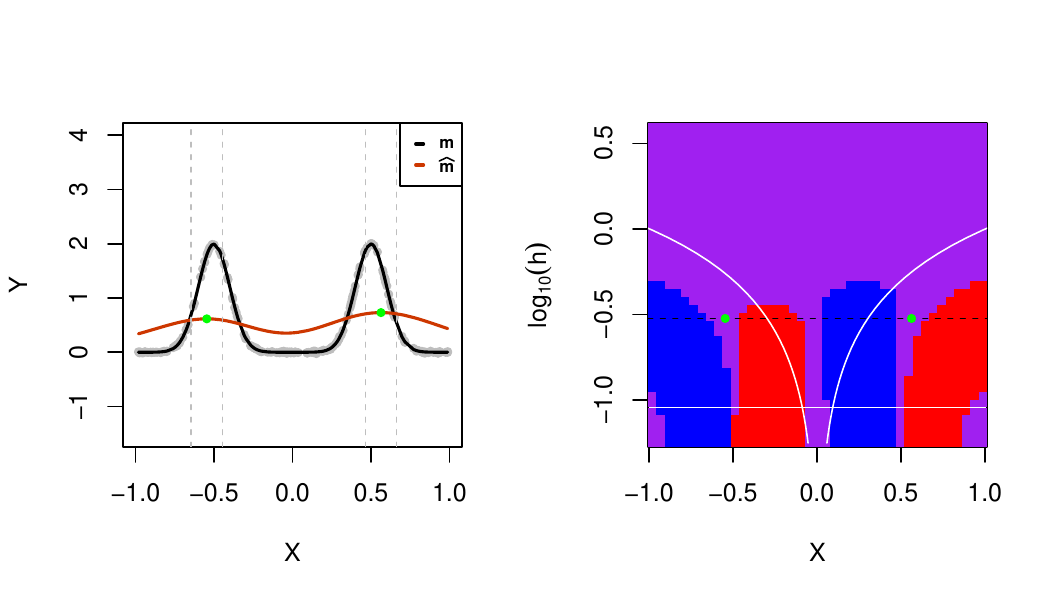}
    \caption{A bimodal function and its estimate, along with the SiZer map of the observed data. Here,  $m(x) = 0.5p_{1}(x) + 0.5 p_{2}(x), x\in [-1, 1]$, where $p_{1}(.)$ and $p_{2}(.)$ are the probability density functions of normal distributions truncated on $[-1, 1]$, each with variance 0.01 and means -0.5 and 0.5, respectively.}
    \label{Fig:Modality}
\end{figure}

We simulate noisy observations from the mixture of two normal distributions truncated on [-1, 1] with means -0.5 and 0.5, with the same variance 0.01, and equal mixture proportions. We choose the kernel bandwidth $h = 0.3$ and estimate the regression function for $n = 200$.   
The graph of the true and estimated function using the truncated Gaussian kernel is given in Figure \ref{Fig:Modality} along with the SiZer map of the data. We obtain critical points of the estimated function numerically, using the bisection method. We focus on two of the critical points depicted by green dots in Figure \ref{Fig:Modality} and aim to infer the modality of the underlying true function in the local regions highlighted by gray vertical dotted lines in Figure \ref{Fig:Modality}. These regions are $ \delta$-neighbourhoods around the critical points of the estimated function, where we choose $\delta = 0.1$. We use the aforementioned characterization to infer whether the critical points represent modes of the true function.

\begin{figure}[htbp!]
    \centering
    \includegraphics[width=\linewidth]{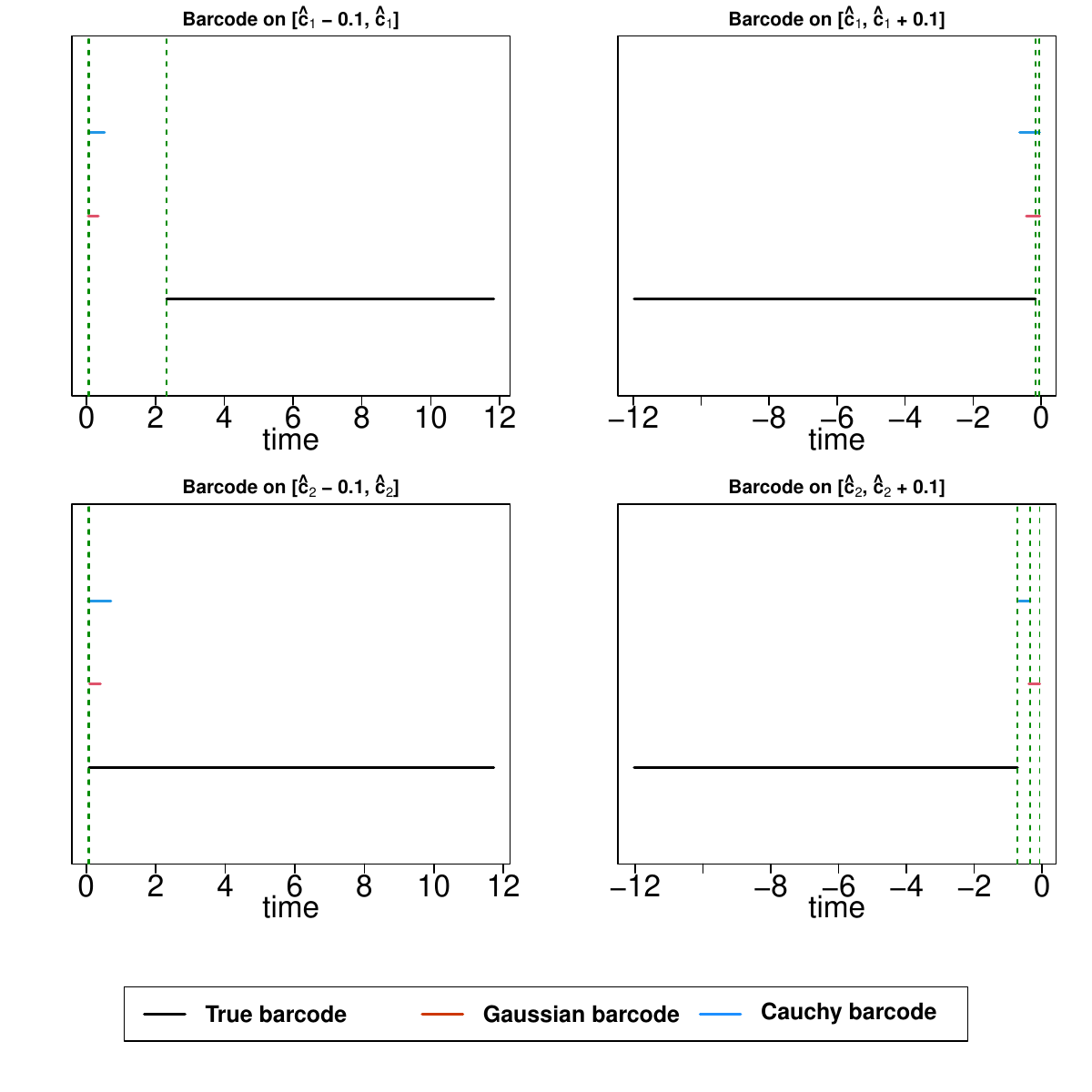},
    \caption{Estimated barcodes of $m_{1}$ obtained from data shown in Figure \ref{Fig:Modality}, in local regions $[\hat{c}_i - 0.1, \hat{c}_i ]$ and $[\hat{c}_i, \hat{c}_i  + 0.1]$, i = 1, 2. Here, $m(x) = 0.5p_{1}(x) + 0.5 p_{2}(x), x\in [-1, 1]$, where $p_{1}(.)$ and $p_{2}(.)$ are the probability density functions of normal distributions truncated on $[-1, 1]$, each with variance 0.01 and means -0.5 and 0.5, respectively.}
    \label{Modality barcodes}
\end{figure}

We obtain barcodes of the estimated derivative using the procedure given in Section \ref{Problem Formulation and Main Results} and infer the monotonicity of the function in the local regions $[c- \delta, c]$ and $[c, c + \delta]$ to infer the modality. Note that here $c$ represents the critical points of interest, and $\delta = 0.1$. We estimate barcodes of $\hat{m}_1$ in these local regions using the estimated critical points. The barcode of $\hat{m}_1$ is obtained using the filtration parameter $\epsilon = 0.005$ for the truncated Gaussian and truncated Cauchy kernel on [-1, 1] using $h = 0.3$ and $n= 200$. We denote the critical points of the estimated function $\hat{m}$ depicted in Figure \ref{Fig:Modality} by $\boldsymbol{\hat{c}_{1}}$ and $\boldsymbol{\hat{c}_{2}}$. The true and estimated barcode are shown in Figure \ref{Modality barcodes}, and the results are articulated in Table \ref{Modality_Results}. 
\begin{table}[htbp!]
\centering
\small
\caption{Results for the truncated Gaussian and Cauchy kernel}
\resizebox{\textwidth}{!}{%
\begin{tabular}{lcc@{\hskip 2pt}cc@{\hskip 6pt}cc}
\toprule
& \multicolumn{2}{c}{\textbf{Estimated critical points ($\hat{c}$)}} 
& \multicolumn{2}{c}{\textbf{Barcode on $[\hat{c} - 0.1, \hat{c}]$}} 
& \multicolumn{2}{c}{\textbf{Barcode on $[\hat{c}, \hat{c} + 0.1]$}} \\
\cmidrule(lr){2-3} \cmidrule(lr){4-5} \cmidrule(lr){6-7}
& \( \mathbf{\hat{c}_{1}} \) & \( \mathbf{\hat{c}_{2}} \) 
& \( [\hat{c}_1 - 0.1, \hat{c}_1] \) & \( [\hat{c}_2 - 0.1, \hat{c}_2] \) 
& \( [\hat{c}_1 , \hat{c}_1 + 0.1] \) & \( [\hat{c}_2 , \hat{c}_2 + 0.1]  \) \\
\midrule
\(\text{Gaussian}\) & -0.54 & 0.56 & \textbf{[0.05, 0.34]} & \textbf{[0.10, 0.41]} & \textbf{[-0.43, -0.04]} & \textbf{[-0.39, -0.07]} \\
\(\text{Cauchy}\)   & -0.51 & 0.52 & \textbf{[0.07, 0.52]} & \textbf{[0.06, 0.71]} & \textbf{[-0.63, -0.05]}  & \textbf{[-0.68, -0.36]} \\
\bottomrule
\end{tabular}%
}
\label{Modality_Results}
\end{table}

Note that we highlight the true and estimated death or birth times in different regions using green dotted lines in Figure~\ref{Modality barcodes} to investigate monotonicity in these regions (see Section~\ref{ Monotonicity}). The observed barcodes in the local regions suggest that the estimated function is non-decreasing in $[\hat{c}_i - 0.1, \hat{c}_i ]$ and non-increasing in $[\hat{c}_i, \hat{c}_i + 0.1 ]$, $i= 1, 2$. This is evident from Table \ref{Modality_Results} and Figure~\ref{Modality barcodes} as the observed death times in $[\hat{c}_i - 0.1, \hat{c}_i ]$ are positive and the observed birth times in $[\hat{c}_i , \hat{c}_i + 0.1 ]$ are negative, $i = 1, 2$. Next, we measure the statistical significance of the observed barcodes by Equation~\eqref{significance measure} using $z_{n, \alpha} = \epsilon$, which is obtained by Equation~\eqref{Znalpha} for $\alpha = 0.05, h = 0.3, \epsilon = 0.005, \text{ and } n  = 200$. Note that here Equation~\eqref{Znalpha} yields $z_{n, \alpha} = \epsilon$ as the constant $C_{\epsilon}$ is large enough, making the estimated barcodes consistent for the aforementioned parameters and sample size. Therefore, as the length of all the observed bars is larger than 0.03, we conclude that the observed bars for both kernels are significant at the $\alpha = 0.05$ level of significance. 

Hence, we conclude that the observed modes are significant at the 95\% level of confidence as per the proposed characterization using barcodes of $\hat{m}_1$. The SiZer map of the data shown in Figure \ref{Fig:Modality} also suggests the same. The estimated critical points (green dots in Figure \ref{Fig:Modality}) lie in a purple region, and the regions on the left and right of these points are colored blue and red, respectively. This suggests the presence of two significant modes in the data for $h = 0.3$, which is highlighted by black dotted lines on the SiZer map in Figure \ref{Fig:Modality}.  

\begin{r1}
    In the context of inferring the modality of a smooth regression function, the following alternative characterization of modality is suggested by one of the reviewers. The reviewer pointed out that one need not consider the estimated derivative $\hat{m}_1$ to characterize the modality; rather, a qualitative comparison of the barcodes of the estimated regression function $\hat{m}$ and -$\hat{m}$ is sufficient. In particular, if $\hat{m}$ has k significant bars and -$\hat{m}$ has at least k + 1 significant bars or -$\hat{m}$ has k significant bars and $\hat{m}$ has at least k + 1 significant bars, conclude that there are at least k modes in the area of interest. This characterization of modality allows us to explore further structures using the derivative of $\hat{m}$ that cannot be revealed by SiZer. For example, given that there is a mode in the region of interest, we can use barcodes of $\hat{m}_1$ and -$\hat{m}_1$ to infer concavity of the function in that region using the proposed characterization of convexity in Section 5.2.    
\end{r1}

\section{Real data study} \label{Real data study}

We implement the proposed methodology on two real data sets to infer local structures such as monotonicity, convexity, and modality. We measure the statistical significance of observed structures using the methodology proposed in Section \ref{Quantitative evaluation}, and compare the results with the respective SiZer maps. Data set 1 is used to implement the monotonicity characterization using the choice of kernel bandwidth obtained from Equation  \eqref{Bandwidth_Selection} to infer the bars at the 5\% level of significance. Data set 2 is used to implement the characterization of convexity and modality using a particular choice of bandwidth and to provide inference using the proposed methodology.

\subsection{Data set 1} \label{Data set 1} 
We consider the cars data available in the R software library, ``datasets.'' This data set consists of $n = 50$ observations on the speed of cars (X) and the stopping distances of cars (Y). We fit a non-parametric regression function to the observed data using the kernel bandwidth obtained from Equation \eqref{Bandwidth_Selection} for $\epsilon = 0.001$ and $\alpha = 0.05$. Note that, to solve Equation \eqref{Bandwidth_Selection} for $h$, we use estimates of $\gamma_1$ (see Assumption B.~\ref{B 5}) and $N_{\epsilon}$ (see Equation \ref{ Theorem constants}) obtained from Algorithm~\ref{Algorithm}. However, in this case $ \eta = Median\{X_{(2)}  - X_{(1)}, \ldots, X_{(n)}  - X_{(n -1)} \} = 0$ in Step 2 of Algorithm~\ref{Algorithm}, which should be > 0. Therefore, we make a slight change in Step 2 (Threshold Selection) of Algorithm~\ref{Algorithm}, and define $ \eta = mean\{X_{(2)}  - X_{(1)}, \ldots, X_{(n)}  - X_{(n -1)} \}$ which is > 0. We refer to Remark \ref{remark:Algo} for more details on the choice of $\eta$ in Algorithm~\ref{Algorithm}. The obtained kernel bandwidth is $h = 15.5$.

\begin{figure}[htbp!]
    \centering
    \includegraphics[width=\linewidth]{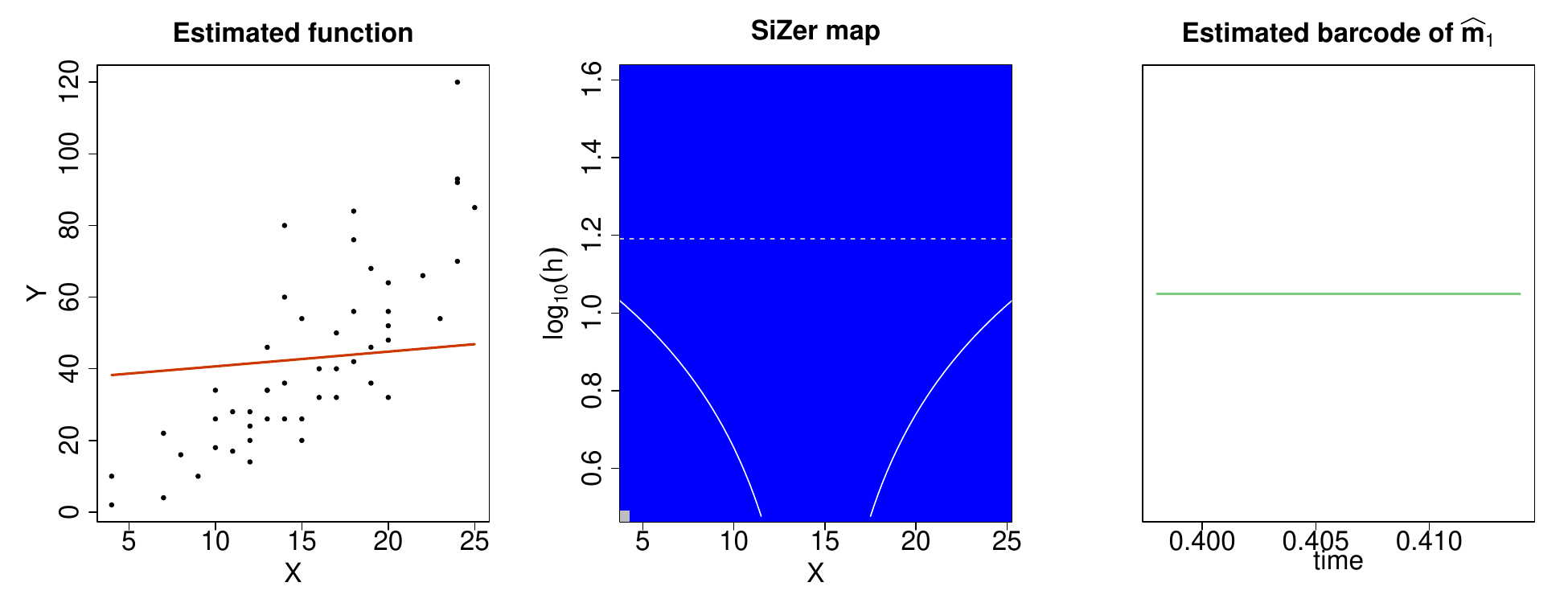}
    \caption{ The observed data (Data set 1), estimated function for $h = 15.5$, SiZer map, and the estimated barcodes of $\hat{m}_1$.}
    \label{cars data plot}
\end{figure}

We fit a non-parametric regression function ($\hat{m}$) to the data, and obtain barcodes of its derivative ($\hat{m}_1$) using the procedure given in Section \ref{Problem Formulation and Main Results}. The estimated barcode of $\hat{m}_1$  is $\{[\boldsymbol{0.398, 0.414}]\}$ which is significant at 5\% level of significance by Equation \eqref{significance measure} for $z_{n, \alpha} = 0.001$. This is because the length of the observed bar is larger than 0.007, and the choice of $z_{n, \alpha}$ comes from the fact that $h = 15.5$ is a solution of Equation \eqref{significance measure} for $\epsilon = 0.001$ and $\alpha = 0.05$. Further, the observed death time is positive, which suggests that the true function underlying the data is strictly increasing. This conclusion is also valid from the SiZer map of the data as shown in Figure \ref{cars data plot}. The gray dotted line on the SiZer map in Figure \ref{cars data plot} highlights $h = 15.5$.          

\subsection{Data set 2}
\begin{figure}[htbp!]
    \centering
    \includegraphics[width=\linewidth]{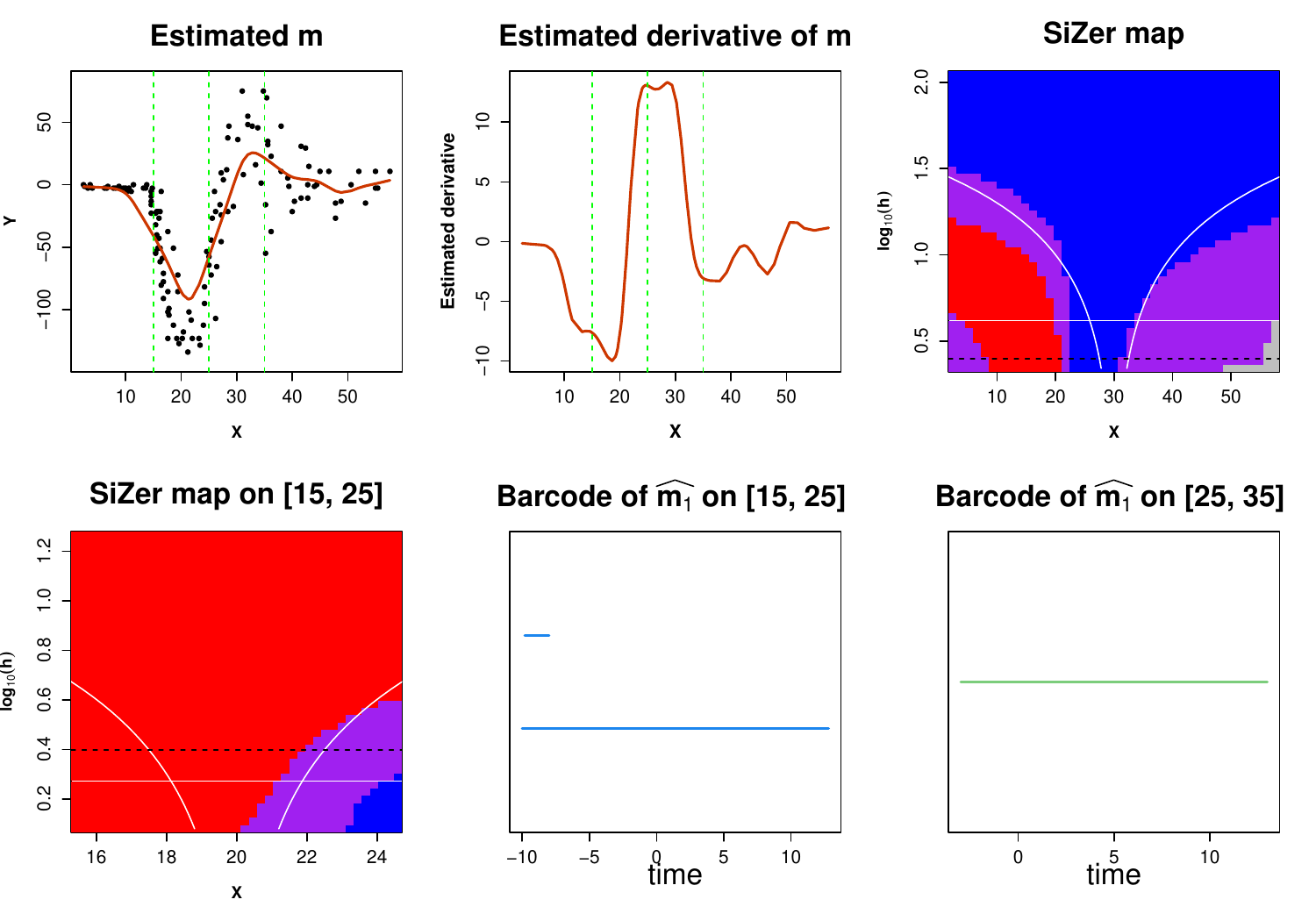}
    \caption{For Data set 2, comparison of local structures in the estimated regression function ($\hat{m}$) derived from the estimated derivative ($\hat{m}_{1}$) with the SiZer map using the bandwidth $h = 2.5$}
    \label{Data set 2}
\end{figure}

We consider the ``Motorcycle data" available in the R package ``adlift". This data set contains $ n = 133$ observations of simulations measuring the effects of motorcycle crashes on victims' heads. There are two variables in the data set: time after a simulated impact with motorcycles (X) and head acceleration of a post-mortem human test object (Y). The observed structure in the data set is shown in Figure \ref{Data set 2}. We fit a non-parametric regression model (see Equation~\ref{Non parametric model}) to the observed data, and aim to infer the local structures present in the data. In particular, we are interested in inferring the observed convexity structure arising in the local region [15, 25], and modality in [25, 35]. These local regions are highlighted by green dotted lines in Figure~\ref{Data set 2}. A non-parametric regression fit ($\hat{m}$), the estimated derivative ($\hat{m}_{1}$), the SiZer map of the data, and the estimated barcodes of $\hat{m}_{1}$ are shown in Figure \ref{Data set 2}. The black dotted line on the SiZer map highlights the kernel bandwidth $h = 2.5$.

We compute the barcode of the estimated derivative using the procedure given in Section \ref{Problem Formulation and Main Results} to investigate whether the aforementioned observed local structures are significant or not. We estimate barcodes for the filtration parameter $\epsilon = 0.1$ and the kernel bandwidth $h = 2.5$. We measure the statistical significance of the observed barcodes using Equation \eqref{significance measure} for $z_{n. \alpha}$ obtained from Equation \eqref{Znalpha} for $\alpha = 0.05$ and for the aforementioned parameters. 

The estimated barcode of $\hat{m}_1$ on [15, 25] is \{\textbf{[-10, 12.8], [-9.8, -8]}\}. Now, by Equation \eqref{Znalpha}, we obtain $z_{n, \alpha} = 0.1$, implying that both observed bars are significant, because the lengths of these bars are larger than 0.7. Therefore, both observed bars are significant at the 95\% level of confidence by Equation \eqref{significance measure}. Hence, by step (i) of the characterization of convexity proposed in Section \ref{Convexity}, we infer that the observed structure in [15, 25] is not convex. This conclusion is also consistent with the following observation pointed out by one of the reviewers. Notice that the two bars in the barcode of the estimated derivative detect a significant valley structure in the estimated derivative in the region [15, 25] (see Figure \ref{Data set 2}). This is because a significant sign change in the estimated derivative is detected by the first bar, and the second bar detects that the estimated derivative has a local minimum in the region [15, 25]. This discards the possibility of the observed structure in the estimated regression function being convex in the region [15, 25].

Moreover, the proposed methodology infers that the observed valley structure in the estimated regression function in the region [15, 25] is statistically significant at the 5\% level of significance. This is because the sign of the estimated derivative is negative at the beginning of the region, following a positive sign in the region [15, 25] (see Figure \ref{Data set 2}). Moreover, this sign change is significant by the first bar in the barcode of the estimated derivative, as the bar is significant and contains 0. However, note that the SiZer map in Figure \ref{Data set 2} does not detect an increasing slope in the region for the chosen bandwidth. Hence, the proposed methodology detects the valley structure underlying the data in the region [15, 25] while the SiZer fails to detect it for the chosen bandwidth.

Next, we estimate the barcode of $\hat{m}_1$ on [25, 35] to investigate the peak structure in the region. The estimated barcode is \{\textbf{[-3, 13]}\}, which is not significant at the 5\% level of significance using the significance measure given in Equation \eqref{significance measure}. Note that there must be a solution of Equation \eqref{Znalpha} for $z_{n, \alpha}$ in $(\epsilon, 16 / 5\sqrt{2})$ for the observed bar to be significant, as the length of the observed bar is 16. Therefore, since there does not exist any solution of Equation \eqref{Znalpha} in $(\epsilon, 16 / 5\sqrt{2})$ for $z_{n, \alpha}$, we conclude that the observed bar is not significant at the 5\% level of significance. Hence, the observed mode in [25, 35] is not significant. The SiZer map in Figure \ref{Data set 2} also suggests the same, as it does not detect the falling slope with confidence in the region [25, 35].            

\section{Conclusion} \label{Conclusion} 

In this article, we explored structures such as monotonicity, convexity, and modality in smooth regression curves via persistent homology of the super-level sets of a particular function. That particular function is the first derivative of the regression function. In general, persistent homology can be applied to identify structures in level sets of functions in terms of connected components and holes. We exploit the fact that derivatives of a function carry important geometric information about the function, and one can identify useful structures in a function by investigating the persistent homology of its derivatives. The main problem that impedes one from using persistent homology to explore the geometric properties of a function via its derivative is the consistency of the estimates of the persistent homology of its derivatives. The problem of consistency for the estimates of persistent homology arises due to the non-availability of an estimator of the derivative of the function that converges uniformly. Therefore, the usual procedure for estimating the persistent homology is ruled out in this case.   

However, \cite{Tdaconsistency} proposed a different procedure based on the kernel estimation of statistical functions (density and regression functions), which is applicable in the present context. We apply this procedure to estimate the persistent homology of the first derivative of regression functions and establish its consistency. This procedure is useful in the present context since its consistency does not require uniform convergence of the estimator of the underlying function. Thus, this procedure paves the way for exploring structures via the persistent homology of the first derivative of the regression functions. For future considerations, it would be tempting to extend this procedure for higher-order derivatives to explore some stimulating structures in regression curves.

We would also like to explore the possibility of extending this procedure for the gradient of the regression function when there is more than one regressor. In this context, we would like to discuss the possibility of extending the procedure when there is more than one regressor and the difficulties involved therein. Note that when there is more than one regressor, the regression function is $m: \mathbb{R}^d \xrightarrow{} \mathbb{R}, d \geq 2$. In this case, the first derivative, that is the gradient of the function, will be a function $m_{1}: \mathbb{R}^d \xrightarrow{} \mathbb{R}^d, d \geq 2 $. The usual persistent homology, also called one-parameter persistent homology, tracks the evolution of homological features in nested sequences of topological spaces indexed by a parameter $t \in \mathcal{I} \subset \mathbb{R}$. For example, filtrations built on point clouds usually have a radius $r \geq 0$ of Euclidean balls, or in super level-set filtrations, usually, it is the height of a real-valued function that varies over $\mathbb{R}$. Therefore, the usual persistent homology can be used only for real-valued functions. In contrast, here we aim to compute the persistent homology of a vector-valued function that does not have the usual notion of height. Thus, to compute persistent homology for the first derivative of a function when the domain of the function is $\mathbb{R}^d, d \geq 2$, we have to resort to the so-called multiparameter persistent homology (see, e.g., \cite{BotnanLesnick2023}). The multiparameter persistent homology is meant to track the evolution of topological features for filtrations of the form $\{\mathcal{X}_{\boldsymbol{t}}; \boldsymbol{t} \in \mathbb{R}^d, \}, d \geq 2$. However, it does not have a straightforward extension of one-parameter persistent homology because it does not have a nice numerical invariant, such as a barcode. Thus, though it provides some information about the underlying topological space, it is not a compact numerical invariant, making it restrictive for topological inference. For these reasons, it is not so straightforward to make topological inferences about the derivative of the regression function when there is more than one regressor.   

\section{Appendix} \label{Appendix}
\subsection{Notations} \label{Notations} 
\begin{itemize}

    \item $a_{n} \asymp b_{n}$ denotes $\lim_{n \to \infty} a_{n}/b_{n} = 1.$
    
    \item$K_{h}(x) = K(\frac{x}{h})$, for any $x \in S$, $S$ is the support of the Kernel function K. 

     \item $K^2(x) \triangleq (K(x))^2 $, for any $x \in S$.

    \item $\left\lceil x \right\rceil$ denotes the smallest integer larger than or equal to $x$.
   
    \item $\mathbb{B}_{h}(x)$ denotes a Ball of radius $h$ centered around the point $x$.

    \item $\mathcal{X}_{n} = \{\left(X_{i}, Y_{i}\right); i = 1,\ldots n\}$

   \item $\partial$A denotes the boundary of set A.

     \item $\partial D_{t}(h) \triangleq \bigcup_{x \in \partial D_{t}} \mathbb{B}_{h}(x)$ where $D_{t}$ denotes super level set at the level t.
     
    \item $D_{t}^{+}(h) \triangleq D_{t} \cup \partial D_{t}(h)$ and $D_{t}^{-}(h) \triangleq D_{t} \setminus\partial D_{t}(h)$.
    
\end{itemize}

\subsection{Definitions} 
\begin{d1} [\textbf{Hoeffding’s inequality}] \label{Hoeffdings}

\noindent Let $Z_{1}, \ldots Z_{n}$ be independent random variables such that $\mathbb{E}(Z_{i}) = 0$, and there exists constants $a_{i} < b_{i}$ such that $\mathbb{P}\left(Z_{i} \in [a_{i}, b_{i}]\right) = 1$ for all $i = 1, \ldots n$. Then for all $u> 0$, we have $$ \mathbb{P}\left(\sum \limits_{i = 1}^{n} Z_{i} \geq u \right) \leq \exp \left( \frac{-2 u^2}{\sum \limits_{i = 1}^{n} \left(b_{i} - a_{i}\right)^2 }\right).$$   
\end{d1}

 \subsection{ Proof of Theorem \ref{Th:1}} \label{Proofs}
    First, we state and prove the following two lemmas that will be used to prove Theorem \ref{Th:1}. Then we provide proof of Theorem \ref{Th:1}. The Lemma \ref{Lemma1} and Lemma \ref{Lemma2} can be interpreted as follows. First, recall the set $\mathcal{X}^{t}_{n} \triangleq \{X_{i} : \hat{f}_{n}(X_{i}) \geq t, i = 1, \ldots, n  \}$ in step 2 of the described procedure. Note that even small errors in the estimator $\hat{f}_{n}$ may result in a false positive or false negative error in the assignments of data points to $\mathcal{X}^{t}_{n}$ with respect to the super-level set that corresponds to the true function $m_{1}$. The true super-level set is denoted as $D_{t}$ here. A tiny error in $\mathcal{X}^{t}_{n}$ may introduce large errors in the estimated persistent homology of $D_{t}$. Therefore, we need to find the probability of false positive and false negative assignments of data points to the set $\mathcal{X}^{t}_{n}$ with respect to $D_{t}$. We find an upper bound of false positive and false negative errors by inflating and deflating the set $D_{t}$ by a radius $h$. In other words, we define the two sets $D_{t}^{+}(h)$ and $D_{t}^{-}(h)$ as follows. First, we define the tube of radius of $h$ around the boundary of $D_{t}$, That is $$\partial D_{t}(h) \triangleq \bigcup_{x \in \partial D_{t}} \mathbb{B}_{h}(x), \partial D_{t} \text{ is boundary of } D_{t} .$$   Then define the following sets: $$ D_{t}^{+}(h) \triangleq D_{t} \cup \partial D_{t}(h) \text{ and } D_{t}^{-}(h) \triangleq D_{t} \setminus\partial D_{t}(h), $$ where $D_{t}^{+}(h)$ and $D_{t}^{-}(h)$ corresponds to inflating and deflating the set $D_{t}$ by a radius $h$, respectively. The Lemma \ref{Lemma1} computes upper bounds for false positive and false negative errors of assignments of data points to the set $\mathcal{X}_{n}^{t}$. Then, using Lemma \ref{Lemma1}, Lemma \ref{Lemma2} shows that the set $\hat{D}_{t} (n, h)$ (See, Equation \ref{hatDnh}) is sandwiched between two non-random approximations of $D_{t}$ with high probability.      

%% Statement of Lemma 1 
\begin{l1} \label{Lemma1}
    Under the assumptions A.~\ref{A 1}-A.~\ref{A 3} and B.~\ref{B 1}-B.~\ref{B 5}, for every $t > 0$ and $\epsilon \in (0, t)$, if $h \equiv h(n) \xrightarrow{} 0 \text{ as } n \xrightarrow{} \infty $ such that $n h^6 \xrightarrow{} \infty$, then there exists a constant $C_{\epsilon}$ such that for large enough $n$ we have
\begin{equation}\label{Eq:3.1}
    \mathbb{P} \left(\exists X_{i} \in \mathcal{X}_{n}^{t}  |  X_{i}\notin D_{t - \epsilon}^{+}(h) \right)\leq n e^{- C_{\epsilon}n h^6}
\end{equation} and 
\begin{equation} \label{Eq:3.2}
    \mathbb{P} \left(\exists X_{i} \notin \mathcal{X}_{n}^{t} |  X_{i}\in D_{t + \epsilon}^{-}(h)\right) \leq n e^{- C_{\epsilon}n h^6},
\end{equation}
where $D_{t}^{+}(h) \triangleq D_{t} \cup \partial D_{t}(h)$, $D_{t}^{-}(h) \triangleq D_{t} \setminus\partial D_{t}(h)$ and $\partial D_{t}(h) \triangleq \bigcup_{x \in \partial D_{t}} \mathbb{B}_{h}(x)$.
\end{l1}

%%%%%%% proof of Lemma 1
\begin{proof} [Proof of Lemma \ref{Lemma1}]
    First, recall that the set $\mathcal{X}_{n}^t = \{ X_{i}: \hat{m}_{1,n}(X_{i}) \geq t\}$, where $\hat{m}_{1,n}$ is defined in Equation \ref{Derivative of m}. Then we start the proof by proving Equation \ref{Eq:3.1}. Consider the probability
\begin{align*}
  \mathbb{P} \left(\exists X_{i} \in \mathcal{X}_{n}^{t}  |  X_{i}\notin D_{t - \epsilon}^{+}(h) \right) & = \bigcup_{i = 1}^{n}\mathbb{P} \left( X_{i} \in \mathcal{X}_{n}^{t}  |  X_{i}\notin D_{t - \epsilon}^{+}(h) \right)\\
&\leq n \mathbb{P} \left( X_{1} \in \mathcal{X}_{n}^{t}  |  X_{1}\notin D_{t - \epsilon}^{+}(h) \right)\\
& = n\mathbb{P} \left( X_{1} \in \mathcal{X}_{n}^{t}  |  X_{1}\in (D_{t - \epsilon}^{+}(h))^c \right)\\
& = n\mathbb{P} \left( \hat{m}_{1, n}(X_{1})\geq t |  X_{1}\in (D_{t - \epsilon}^{+}(h))^c \right)\\
& = n\int_{(D_{t - \epsilon}^{+}(h))^c} \int_{\mathbb{R}}f(x,y)\mathbb{P}(\hat{m}_{1, n}(X_{1}) \geq t| X_{1} = x, Y_{1} = y) dy dx.\tag{E.2}\label{Eq:E.2}    
\end{align*}
 %where $$\hat{m}_{1, n}(x) = \frac{\mathrm{d}}{\mathrm{d}x}\hat{m}_{n}(x) = \frac{\sum\limits_{i=1}^n K_{h}\left(x-X_{i}\right) \sum\limits_{i=1}^n Y_{i}K_{1,h}\left(x-X_{i}\right) - \sum\limits_{i=1}^n Y_{i}K_{h}\left(x-X_{i}\right)\sum\limits_{i=1}^n K_{1,h}\left(x-X_{i}\right) }{\left(\sum\limits_{i=1}^n K_{h}(x-X_{i})\right)^2},$$  
%$$ \hat{m}_{n}(x) = \frac{\sum\limits_{i = 1}^{n} K_{h}(x - X_{i}) Y_{i}}{\sum\limits_{i = 1}^{n} K_{h}(x - X_{i})}, K_{1,h}(x-X_{i}) = \frac{\mathrm{d}}{\mathrm{d}x} K_{h}(x-X_{i})\text{, } K_{h}(x) = K(x/h).$$
Now, consider the following:
\begin{align*}
    \hat{m}_{1, n}(x) &=  \frac{\sum\limits_{i=1}^n K_{1,h}(x-X_{i}) \left( Y_{i} - \hat{m}_n(x)\right)}{\sum\limits_{i=1}^n K_{h}(x-X_{i})}\\
    &\leq \frac{\sum\limits_{i=1}^n K_{1,h}(x-X_{i}) \left( Y_{i} - Y_{min}\delta\right)}{\sum\limits_{i=1}^n K_{h}(x-X_{i})} \text{ almost surely }, \tag{E.3} \label{mhatn}
\end{align*}
where \ref{mhatn} follows from the fact that $\hat{m}_{n}(x) \geq Y_{min}\delta \text{ almost surely for all } x \in \mathbb{R}$. 

Now, consider the probability
\begin{align*}
    \mathbb{P}\left[\hat{m}_{1, n}(X_{1}) \geq t| X_{1} = x, Y_{1} = y\right]
    & = \mathbb{P}\left[\frac{\sum\limits_{i=1}^n K_{1,h}\left(X_{1}-X_{i}\right) \left(Y_{i} - \hat{m}_{n}(X_{1})\right)}{\sum\limits_{i=1}^n K_{h}\left(X_{1}-X_{i}\right)} \geq t|X_{1} = x, Y_{1} = y\right  ]\\ 
    &\leq \mathbb{P}\left[\frac{\sum\limits_{i=1}^n K_{1,h}\left(X_{1}-X_{i}\right) \left(Y_{i} - Y_{min} \delta)\right)}{\sum\limits_{i=1}^n K_{h}\left(X_{1}-X_{i}\right)} \geq t|X_{1} = x, Y_{1} = y\right  ] \tag{Using \ref{mhatn}}\\  
    & = \mathbb{P}\left[\sum\limits_{i = 2}^{n} \left\{K_{1,h}\left(x-X_{i}\right) \left( Y_{i} -Y_{min} \delta\right) - tK_{h}\left(x-X_{i}\right)\right\} \geq t\right]\\
    &\leq \mathbb{P}\left[\sum \limits_{i = 2}^{n} \left(Z_{i} - \mathbb{E}\left(Z_{i}\right)\right) \geq t - \sum \limits_{i = 2}^{n}\mathbb{E}\left(Z_{i}\right)\right] \\
    &\leq \mathbb{P}\left[\sum \limits_{i = 2}^{n} \left(Z_{i} - \mathbb{E}\left(Z_{i}\right)\right) \geq t + 2h\delta p_{min} ( \gamma_1 + \epsilon) (n - 1)\right], \tag{E.4 } \label{Eq:E.4}
     \end{align*}
where \begin{equation} \label{Zi}
    Z_{i} = K_{1,h}\left(x-X_{i}\right)\left(Y_{i} - Y_{min} \delta \right) - (m_{1}(X_{i}) + \epsilon) K_{h}\left(x-X_{i}\right).
\end{equation} We have defined $Z_{i}$ using the fact that since $X_{i} \notin D^{+}_{t - \epsilon}(h)$ this implies that $m_{1}(X_{i}) \leq t - \epsilon$, and the inequality in Equation \ref{Eq:E.4} follows from the Result \ref{Result 1}.  

Now, we would like to apply Hoeffding's inequality (See Definition \ref{Hoeffdings}) for the random variables $Z_{i} - \mathbb{E}Z_{i}, i = 2,\ldots n$, therefore consider the following:
\begin{enumerate}
       \item We have $\mathbb{E}[Z_{i} - \mathbb{E}Z_{i}] = 0$, and $Z_{i} - \mathbb{E}Z_{i}$'s are independent, since $(X_{i}, Y_{i})$'s are independent for all $i = 2, \ldots n$.

    \item We have $Z_{min}\leq Z_{i} \leq Z_{max}$, almost surely, where $Z_{min}$ and $Z_{max}$ are defined as:
$$ Z_{min} = (X_{min} - X_{max}) h^{-2} \tau_1 Y_{min} \delta (1 - \delta ) - (M_{1} + \epsilon),$$
$$ Z_{max} = (X_{max} - X_{min}) h^{-2} \tau_2 (Y_{max} - Y_{min}\delta ) - (\gamma_{1} + \epsilon) \delta.$$
This follows from Assumptions A.~\ref{A 1}, A.~\ref{A 3}, B.~\ref{B 1}, B.~\ref{B 3}, and B.~\ref{B 5}. Consequently, assuming $ (M_{1} + \epsilon)h^2 \approx 0$ and $(\gamma_{1} + \epsilon) \delta h^2 \approx 0$, for small enough $h$, we have:
\begin{equation} \label{Zbounds}
    Z_{max} - Z_{min} \approx h^{-2} Q, 
\end{equation}
where $Q = (X_{max} - X_{min}) \left[\tau_2 (Y_{max} - Y_{min}\delta ) + \tau_1 Y_{min} \delta (1 - \delta ) \right]$. 
    \end{enumerate}

Therefore, using 1 and 2 above, and applying Hoeffding's inequality for $ u = t + 2h\delta p_{min} ( \gamma_1 + \epsilon)(n - 1)$, we have
\begin{align*}
 \mathbb{P}\left[\hat{m}_{1, n}(X_{1}) \geq t| X_{1} = x, Y_{1} = y\right] 
 &\leq \exp{ \frac{-2 \left( t + 2h\delta p_{min} ( \gamma_1 + \epsilon)(n - 1)\right)^2 }{\left(n -1\right)\left(Z_{max} - Z_{min}\right)^2} }\\
 & \approx \exp{ \frac{-2n^2 h^2 \left(\frac{t}{nh} + 2\delta p_{min} ( \gamma_1 + \epsilon) \frac{(n - 1)}{n} \right)^2 }{(n-1) \left(h^{-2} Q\right)^2} } \tag{By Equation \eqref{Zbounds}} \\
 & \asymp \exp{ \frac{-2nh^2 (2\delta p_{min} ( \gamma_1 + \epsilon))^2 }{ h^{-4} Q^2} } \\
 & = \exp{(-nh^6 C_{\epsilon})} \tag{E.5}\label{Exponential bound}, 
\end{align*}
 where $C_{\epsilon} = \displaystyle{8\delta^2 p_{min}^2 ( \gamma_1 + \epsilon)^2 Q^{-2}, \text{ provided }  \gamma_1 \neq -\epsilon} $.

 Now, using Equation \ref{Exponential bound} in Equation \ref{Eq:E.2} we have 
 \begin{align*}
  \mathbb{P} \left(\exists X_{i} \in \mathcal{X}_{n}^{t}  |  X_{i}\notin D_{t - \epsilon}^{+}(h) \right) 
  & \leq n\int_{(D_{t - \epsilon}^{+}(h))^c} \int_{\mathbb{R}}f(x,y) \exp{(-nh^6 C_{\epsilon})}dy dx\\
  & = n \exp{(-nh^6 C_{\epsilon})}\int_{\mathbb{R}} \int_{\mathbb{R}}f(x,y)dy dx\\
  & = n e^{-nh^6 C_{\epsilon}}
 \end{align*}
 This completes the proof of inequality in Equation \ref{Eq:3.1}. To prove the inequality in \ref{Eq:3.2}, we proceed in the same way using the fact that $X_{i} \in D_{t + \epsilon}^{-}(h)$ implies that $m_{1}(X_{i}) \geq t + \epsilon$, and define random variables $Z_{i}^*$ as $$ Z_{i}^* = \left(m_{1}(X_{i}) - \epsilon\right) K_{h}(x - X_{i}) - K_{1,h}(x - X_{i}) \left(Y_{i} - Y_{min} \delta\right).$$
 Now, observe that similar to Result \ref{Result 1}, we will have same result for $\mathbb{E}[Z_{i}^*]$, thus similar to Equation \ref{Eq:E.4}, we have $$\mathbb{P}\left[\hat{m}_{1, n}(X_{1}) \leq t| X_{1} = x, Y_{1} = y\right] \leq \mathbb{P}\left[\sum \limits_{i = 2}^{n} \left(Z_{i}^* - \mathbb{E}\left(Z_{i}^*\right)\right) \geq 2h\delta p_{min} ( \gamma_1 + \epsilon) (n - 1) - t\right] .$$ Then apply Hoeffding's inequality for random variables $Z_{i}^* - \mathbb{E}[Z_{i}^*]$ and proceed in the same way to get Equation \ref{Exponential bound} and, from \ref{Exponential bound}, it is easy to prove the inequality in Equation \ref{Eq:3.2}. 
 
 This completes the proof of Lemma \ref{Lemma1}.
\end{proof}

%%% Lemma 2 
\begin{l1} \label{Lemma2}
  Under the assumptions A.~\ref{A 1}-A.~\ref{A 3} and B.~\ref{B 1}-B.~\ref{B 5}, for every $t > 0$, and $\epsilon \in (0, t)$ if $h \equiv h(n) \xrightarrow{} 0 \text{ as } n \xrightarrow{} \infty \text{ such that } nh^6 \xrightarrow{} \infty$, then for large enough $n$, we have 
\begin{equation} \label{Eq:Lemma2}
    \mathbb{P}(D^{-}_{t + \epsilon}(2h) \subset \hat{D}_{t}(n, h) \subset D^{+}_{t - \epsilon}(2h)) \geq 1 - 3n e^{-C_{\epsilon}n h^6}.
\end{equation}  
\end{l1}

\begin{proof} [Proof of Lemma \ref{Lemma2}]
    Let $E_{1}$ and $E_{2}$ denotes the following events respectively,  $$D^{-}_{t + \epsilon}(2h) \subset \hat{D}_{t}(n, h) \text{ and } \hat{D}_{t}(n, h) \subset D^{+}_{t - \epsilon}(2h)).$$
Then, the probability in Equation \ref{Eq:Lemma2} can be written as: 
\begin{equation} \label{Eq:3.6}
\mathbb{P}\left(E_{1} \cap E_{2} \right) = 1 -
 \mathbb{P}\left(E_{1}^c \cup E_{2}^c \right)\geq 1 - \mathbb{P}\left(E_{1}^c\right) - \mathbb{P}\left(E_{2}^c\right).
\end{equation}
Now, since $\hat{D}_{t}(n, h) \supset \mathcal{X}_{n}^t$, using Lemma \ref{Lemma1} we have 
\begin{align*}
    \mathbb{P}\left(E_{2}^c\right) &=\mathbb{P}\left(\hat{D}_{t}(n, h) \not\subset D^{+}_{t - \epsilon}(2h))\right)
    =\mathbb{P}\left(\hat{D}_{t}(n, h)\subset \left(D^{+}_{t - \epsilon}(2h)\right)^c\right)\\
    &\leq \mathbb{P} \left(\mathcal{X}_{n}^t \cap \left(D^{+}_{t - \epsilon}(h)\right)^c \neq \varnothing \right) \leq n e^{-nh^6C_{\epsilon}} . \tag{E.6} \label{First bound}
\end{align*}
 
 Now, consider the probability,
\begin{equation}\label{E1complement}
    \mathbb{P}\left(E_{1}^c\right) = \mathbb{P}\left(D^{-}_{t + \epsilon}(2h)\not\subset \hat{D}_{t}(n, h)\right)= \mathbb{P}\left(\exists S \subset D^{-}_{t + \epsilon}(2h): S \cap \hat{D}_{t}(n, h) = \varnothing \right).
\end{equation}
Note that in view of the fact that $m_{1}$ is bounded, $D^{-}_{t + \epsilon}(2h)$ is bounded, hence it can be covered by finitely many subsets of it. Let $\gamma \in \left(0, 1\right)$ and define a cover of $D^{-}_{t + \epsilon}(2h)$ as  $\left\{ \mathbb{B}_{\gamma h}(u) : u \in \mathcal{I}\right\}$, where $\mathcal{I} \subset D^{-}_{t + \epsilon}(2h) $ is a finite index set. Then for all $x \in D^{-}_{t + \epsilon}(2h)$ there exists a $u \in \mathcal{I}$ such that $|x - u|\leq \gamma h$. By the construction of $\hat{D}_{t}(n, h)$, if $x \not\in \hat{D}_{t}(n, h)$ then $x \not \in \mathbb{B}_{h}(X_{i})$ for all $X_{i} \in \mathcal{X}_{n}^{t}$.

\noindent Therefore, for any $ x \in S \subset D^{-}_{t + \epsilon}(2h)$, where S satisfies that $S \cap \hat{D}_{t}(n, h) = \varnothing$, there exists a $u \in \mathcal{I}$ such that $|x - u| > (1- \gamma) h$. This implies that for any $x \in D^{-}_{t + \epsilon}(2h)$ such that $x \not \in \hat{D}_{t}(n, h)$ there exists a $ u \in \mathcal{I}$ such that $\mathbb{B}_{(1 - \gamma)h}(u) \cap \hat{D}_{t}(n, h) = \varnothing.$

Thus, the probability in Equation \ref{E1complement} can be written as 
\begin{align*}
    \mathbb{P}\left(E_{1}^c\right) &= \mathbb{P}\left(\exists u \in \mathcal{I} : \mathbb{B}_{(1 - \gamma)h}(u) \cap \hat{D}_{t}(n, h) = \varnothing \right)\\
    &\leq \mathbb{P}\left(\exists u \in \mathcal{I} : \mathbb{B}_{(1 - \gamma)h}(u) \cap \mathcal{X}^{t}_{n} = \varnothing \right), \tag{E.7}\label{E.7}
\end{align*}
where Equation \ref{E.7} follows from the fact that $\mathcal{X}^{t}_{n} \subset \hat{D}_{t}(n, h) $.   
Now, note that there could be points in $\mathcal{X}^{t}_{n}$, that are added due to the error in the estimator $\hat{m}_{1,n}$, therefore probability in Equation \ref{E.7} can be written as 
\begin{align*}
 \mathbb{P}\left(E_{1}^c\right) &\leq \mathbb{P}\left(\exists u \in \mathcal{I} : \mathbb{B}_{(1 - \gamma)h}(u) \cap \mathcal{X}^{t}_{n} = \varnothing \right) \\
 &= \mathbb{P}\left(\exists u \in \mathcal{I} : \mathbb{B}_{(1 - \gamma)h}(u) \cap \mathcal{X}^{t}_{n} = \varnothing; D^{-}_{t + \epsilon}(h) \cap \mathcal{X}_{n} \subset \mathcal{X}^{t}_{n} \right)\\
 & + \mathbb{P}\left(\exists u \in \mathcal{I} : \mathbb{B}_{(1 - \gamma)h}(u) \cap \mathcal{X}^{t}_{n} = \varnothing; D^{-}_{t + \epsilon}(h) \cap \mathcal{X}_{n} \not \subset \mathcal{X}^{t}_{n} \right)\\
 & \leq \mathbb{P}\left(\exists u \in \mathcal{I} : \mathbb{B}_{(1 - \gamma)h}(u) \cap \mathcal{X}^{t}_{n} = \varnothing\right) + \mathbb{P}\left(D^{-}_{t + \epsilon}(h) \cap \mathcal{X}_{n} \not \subset \mathcal{X}^{t}_{n} \right). \tag{E.8} \label{E.8}
\end{align*}
Now, consider the following:
\begin{align*}
    \mathbb{P}\left(\exists u \in \mathcal{I} : \mathbb{B}_{(1 - \gamma)h}(u) \cap \mathcal{X}^{t}_{n} = \varnothing\right) &= \sum \limits_{u \in \mathcal{I}} \mathbb{P}\left(\mathbb{B}_{(1 - \gamma)h}(u) \cap \mathcal{X}_{n} = \varnothing\right) \\
    & = \sum \limits_{u \in \mathcal{I}} \mathbb{P}\left(X_{i} \not \in \mathbb{B}_{(1 - \gamma)h}(u) \text{ for all i }\in \mathcal{X}_{n} \right)\\
    & = \sum \limits_{u \in \mathcal{I}} \left(1 - F_{X}\left(\mathbb{B}_{(1 - \gamma)h}(u)\right) \right)^n\\
    &\leq \sum \limits_{u \in \mathcal{I}} e^{-nF_{X}\left(\mathbb{B}_{(1 - \gamma)h}(u)\right)}, \tag{E.9} \label{E.9}
\end{align*}
where $F_{X}$ is the distribution function of the regressors. Then, for $0 < h < 1$, we have 
\begin{equation} \label{Bound for E.9}
    F_{X}\left(\mathbb{B}_{(1 - \gamma)h}(u)\right) = \int\limits_{u - (1 - \gamma)h}^{u +(1 - \gamma)h} p(x) dx \geq 2(1 - \gamma)h p_{min} \geq 2(1 - \gamma)h^6 p_{min} 
\end{equation}
Also, note that since $\mathcal{I}$ is a finite set, there exist finitely many intervals of length $2\gamma h$ that cover the set $D^{-}_{t + \epsilon}(h)$. This implies that there exists a positive number $N_{\gamma}$ such that $|\mathcal{I}| \leq N_{\gamma} 2\gamma h$. Now using this fact with Equation \ref{Bound for E.9} in Equation \ref{E.9}, we have
\begin{equation}  \label{First term bound} 
    \mathbb{P}\left(\exists u \in \mathcal{I} : \mathbb{B}_{(1 - \gamma)h}(u) \cap \mathcal{X}^{t}_{n} = \varnothing\right) \leq  N_{\gamma} 2\gamma h e^{-n 2(1 - \gamma)h^6 p_{min}} = N_{\gamma} 2\gamma h e^{-n h^6 C_{\gamma}},
\end{equation}
where $ C_{\gamma} = 2(1 - \gamma) p_{min} $. Now, note that the second term in Equation \ref{E.8} can be bounded by Lemma \ref{Lemma1}. Therefore, using Equation \ref{Eq:3.2}, we have the following bound for the second term of Equation \ref{E.8}, \begin{equation} \label{Second term bound}
    \mathbb{P}\left(D^{-}_{t + \epsilon}(h) \cap \mathcal{X}_{n} \not \subset \mathcal{X}^{t}_{n} \right) \leq n e^{- C_{\epsilon}n h^6}
\end{equation}     
Now, using Equation \ref{First term bound} and Equation \ref{Second term bound} in Equation \ref{E.8}, we have 
\begin{equation} \label{ Final bound}
    \mathbb{P}\left(E_{1}^c\right) \leq N_{\gamma} 2\gamma h e^{-n h ^2 C_{\gamma}} + n e^{- C_{\epsilon}n h^2} \leq 2n e^{- C_{\epsilon}n h^6},
\end{equation}
where the last inequality follows from the observation that $C_{\epsilon} < C_{\gamma}$, and for large n, $ N_{\gamma} 2\gamma h \leq n $.

Thus, using Equation \ref{ Final bound} and Equation \ref{First bound} in Equation \ref{Eq:3.6}, we get the desired inequality as in Equation \ref{Eq:Lemma2}.
This completes the proof of Lemma \ref{Lemma2}.

\end{proof}

%%%%%%%%%%%%%%% proof of Theorem 3.2

\begin{proof} [Proof of Theorem \ref{Th:1}]
    Let $PH_{k}^\epsilon(m_{1})$ denotes the persistent homology of the filtration $\mathcal{D}^\epsilon$ defined as 
$$\mathcal{D}^\epsilon \triangleq \{D_{t_{i} + \epsilon}\}_{i\in \mathbb{Z}}, \text{ where } t_{max} = 2\epsilon N_{\epsilon}, t_{i} = t_{max} - 2 i \epsilon \text{, and } N_{\epsilon} =  \sup_{x \in \mathbb{R}} \left\lceil m_{1}(x)/2\epsilon \right\rceil.$$
Since $\mathcal{D}^{\epsilon}$ is a discrete approximation of the continuous filtration $\mathcal{D}$ with step size 2$\epsilon$, therefore we have the maximum difference between $PH_{k}(m_{1})$ and $PH_{k}^\epsilon(m_{1})$ to be 2$\epsilon$. That is, we have  
\begin{equation} \label{Eq:1}
\delta_{B}\left(PH_{k}(m_{1}), PH_{k}^\epsilon(m_{1}) \right) \leq 2\epsilon
\end{equation}
Therefore, it is enough to prove that $\delta_{B}\left(\widehat {PH}_{k}^\epsilon (m_{1}), PH_{k}^\epsilon(m_{1}) \right) \leq 3\epsilon$ holds with high probability. Let E be the event that the following sequence of inclusions holds for $i = 1, \ldots, M_\epsilon$: 
$$D_{t_{i} + \epsilon} \hookrightarrow \hat{D}_{t_{i}}(n, h) \hookrightarrow D_{t_{(i+ 1)} + \epsilon},$$
where $M_\epsilon = \max\{i \in \mathbb{N} : t_i \geq \gamma_1  \}$. Recall that $\gamma_1 = \inf m_1(x)$, defined in Assumption B.~\ref{B 5}. Note that for all $i \geq M_\epsilon$, the level sets $D_{t_{i}}$ remain the same. 

Now, applying Lemma \ref{Lemma2} $M_{\epsilon}$ times for $\hat{D}_{t_{i}}(n, h)$, we have 
$$\mathbb{P}(E) \geq 1 - 3n M_{\epsilon} e^{ -C_{\epsilon/2}n h^6}$$
Therefore, $\mathcal{D}^{\epsilon}$ and $\hat{\mathcal{D}}^{\epsilon}$ are weakly $\epsilon$- interleaving. (See Definition 4.1 \cite{Proximityofpd}). Now, using the weak stability theorem (\cite{Proximityofpd}), we have 
\begin{equation} \label{Eq:2}
    \delta_{B}\left(\widehat{PH}_{k}^{\epsilon}(m_{1}), PH_{k}^{\epsilon}(m_{1}) \right) \leq 3\epsilon.
\end{equation}
Thus, from Equation \ref{Eq:1} and Equation \ref{Eq:2} we have $$\delta_{B}\left(\widehat{PH}_{k}^{\epsilon}(m_{1}), PH_{k}(m_{1}) \right) \leq 5\epsilon.$$

Let $E_{1}$ denotes the event that $\delta_{B}\left(\widehat{PH}_{k}^{\epsilon}(m_{1}), PH_{k}(m_{1}) \right) \leq 5\epsilon$ then since $E\subset E_{1}$ we have $$\mathbb{P}\left( \delta_{B} \left(\widehat{PH}_{k}^{\epsilon}\left(m_{1}\right), PH_{k}\left(m_{1}\right)\right) \leq 5\epsilon\right) \geq 1 - 3nM_{\epsilon}e^{-C_{\epsilon/2}n h^6}.$$
This completes the proof of Theorem \ref{Th:1}.
\end{proof}

\subsection{Some Auxiliary Results} 
%% Result 1
\begin{a1} \label{Result 1}
 For the random variable $Z_{i}$ defined in Equation \ref{Zi}, we have: 
$$\mathbb{E}[Z_{i}] \leq -2h\delta p_{min} ( \gamma_1 + \epsilon)$$.   
\end{a1} 

%% Proof of the Result 1 
\begin{proof}
    Consider the following:
\begin{align*}
    \mathbb{E}[Z_{i}] &= \mathbb{E}_{X_{i}} \left[K_{1,h}\left(x - X_{i}\right)\mathbb{E}_{Y_{i}|X_{i}}\left(Y_{i} - m\left(x\right)\right)\right] - \mathbb{E}_{X_{i}} \left[\left(m_{1} \left(X_{i} \right)+ \epsilon \right)K_{h}(x - X_{i})\right]\\
   &= \int_{x -  h}^{x +  h} K_{1,h}(x - \xi) \left(m(\xi) - m(x)\right)p(\xi) d\xi - \int_{x-h}^{x+h} \left(m_{1} \left(\xi \right)+ \epsilon \right)K_{h}(x - \xi) p(\xi) d\xi \\
   &\leq p_{max} \left(M - m(x)\right) \int_{x- h}^{x +  h} K_{1,h}(x - \xi)d\xi - \delta p_{min}\int_{x - h}^{x +h} m_{1}\left(\xi\right) d\xi - \epsilon \delta p_{min} 2h  \\ 
   &= p_{max} \left(M - m(x)\right)\left[K_{h}\left( - h\right) - K_{h}\left( h\right) \right] -  \delta p_{min}\left[m(x + h) - m( x -h)\right] - \epsilon \delta p_{min} 2h\\
   & = - \delta p_{min}\left[m(x + h) - m( x -h)\right] - \epsilon \delta p_{min} 2h \tag{By Assumption A.~\ref{A 1}} \\
   &\leq -2h\delta p_{min} ( \gamma_1 + \epsilon),
\end{align*}
where the last line follows from the application of the mean value theorem for $m$ on the interval $(x - h, x + h)$, and then subsequently using Assumption B.~\ref{B 5}.
\end{proof}

 \subsection{Comparison and Implementation of Algorithm~\ref{Algorithm}} \label{Alogorithm 1}
    In this subsection, we implement Algorithm~\ref{Algorithm} on the simulated data sets used in Section \ref{Simulation studies} to investigate whether the estimated derivative values are reasonable or not. Note that, for our purpose, some reasonable bounds of the unknown population quantities are sufficient rather than the exact values. We compare the derivative values estimated from Algorithm~\ref{Algorithm} with the true values and values estimated from the ratio of first-order differences of Y(response) and X(regressor). We denote $\hat{Y}_1^T, \hat{Y}_1^A \text{ and } \hat{Y}_1^N$ as the sets that contain true derivative values, derivative values estimated from Algorithm~\ref{Algorithm}, and derivative values estimated from the ratio of first-order differences of Y and X, respectively. 

We simulate $n = 1000$ observations from the functions used in Section~\ref{ Monotonicity}, Section~\ref{Convexity}, and Section~\ref{Modality} with the same error distribution. We denote simulated data sets, for $n = 1000$, used in Section \ref{ Monotonicity}, Section \ref{Convexity}, and Section \ref{Modality} as $D_1$, $D_2 \text{, and } D_3$, respectively. We provide summary statistics of the estimated derivative values and plot the estimated values for visual comparison with the true values. We denote the first and third quartiles by $\boldsymbol{Q_1}$ and $\boldsymbol{Q_3}$, respectively, in the summary tables. The following are the results for the data sets $D_1, D_2 \text{ and } D_3$.   

\begin{figure}[htbp!]
    \centering
    \includegraphics[width=\linewidth]{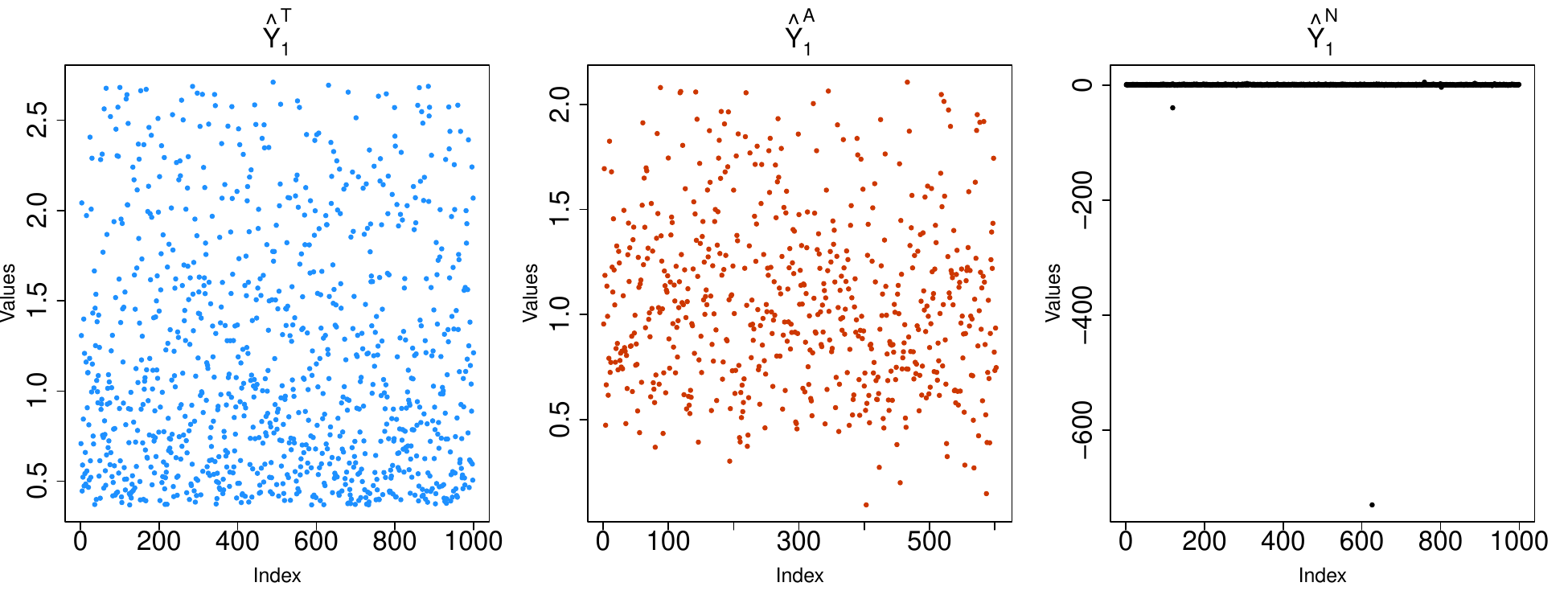}
    \caption{ Plots for $D_1$. }
    \label{D_1}
\end{figure}

\textbf{Results for $\boldsymbol{D_1}$}: 
It is evident from Table \ref{Table:D1} and Figure \ref{D_1} that Algorithm~\ref{Algorithm} provides good estimates except for the lower extreme of the true values. On the other hand, the values in $\hat{Y}_1^N$ are too noisy.     

\begin{table}[htbp!] 
\centering
\caption{Summary statistics for $D_1$}
\label{Table:D1}
\begin{tabular}{lcccccc}
\toprule
\textbf{Values} & \textbf{Min.} & $\boldsymbol{Q_1}$ & \textbf{Median} & \textbf{Mean} & $\boldsymbol{Q_3}$ & \textbf{Max.} \\
\midrule
$\hat{Y}_1^T$    & 0.3693 & 0.6018 & 0.9914 & 1.1620 & 1.6033 & 2.7108 \\
$\hat{Y}_1^A$  & 0.0945 & 0.7538 & 1.0169 & 1.0621 & 1.3013 & 2.1057 \\
$\hat{Y}_1^N$    & -730.0287 & 0.7685 & 1.0332 & 0.3342 & 1.3647 & 5.6000 \\
\bottomrule
\end{tabular}
\end{table}

\begin{figure}[htbp!]
    \centering
    \includegraphics[width=\linewidth]{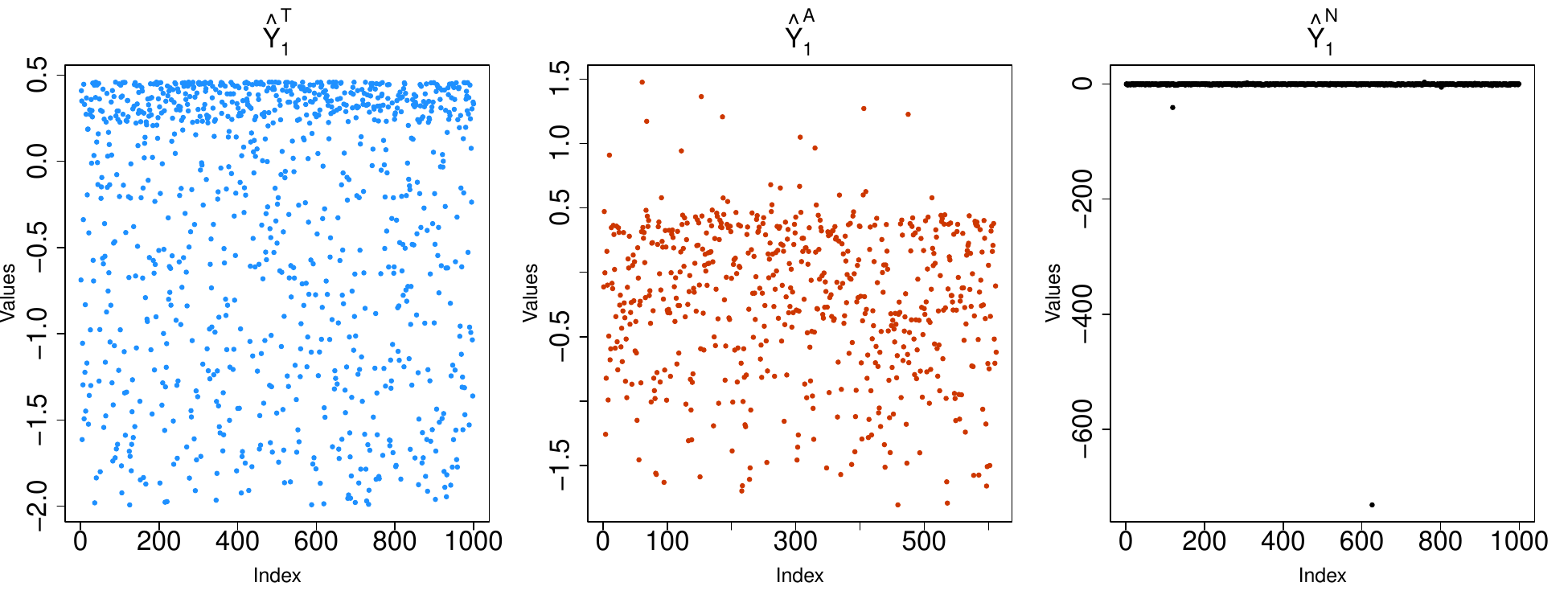}
    \caption{ Plots for $D_2$. }
    \label{D_2}
\end{figure}

\textbf{Results for $\boldsymbol{D_2}$}: 
It is evident from Table \ref{Table:D2} and Figure \ref{D_2} that Algorithm~\ref{Algorithm} provides good estimates except for the upper extreme of the true values. On the other hand, the values in $\hat{Y}_1^N$ are too noisy in this case as well.          
\begin{table}[htbp!]
\centering
\caption{Summary statistics for $D_2$}
\label{Table:D2}
\begin{tabular}{lcccccc}
\toprule
\textbf{Values} & \textbf{Min.} & $\boldsymbol{Q_{1}}$ & \textbf{Median} & \textbf{Mean} & $\boldsymbol{Q_3}$ & \textbf{Max.} \\
\midrule
$\hat{Y}_1^T$    & -1.9925 & -1.0156 & -0.0173 & -0.3357 & 0.3503 & 0.4593 \\
$\hat{Y}_1^A$    & -1.8059 & -0.5828 & -0.1350 & -0.2103 & 0.2641 & 1.4762 \\
$\hat{Y}_1^N$    & -731.5696 & -0.5851 & -0.1371 & -1.0017 & 0.2696 & 3.3662 \\
\bottomrule
\end{tabular}
\end{table}

\begin{figure}[htbp!]
    \centering
    \includegraphics[width=\linewidth]{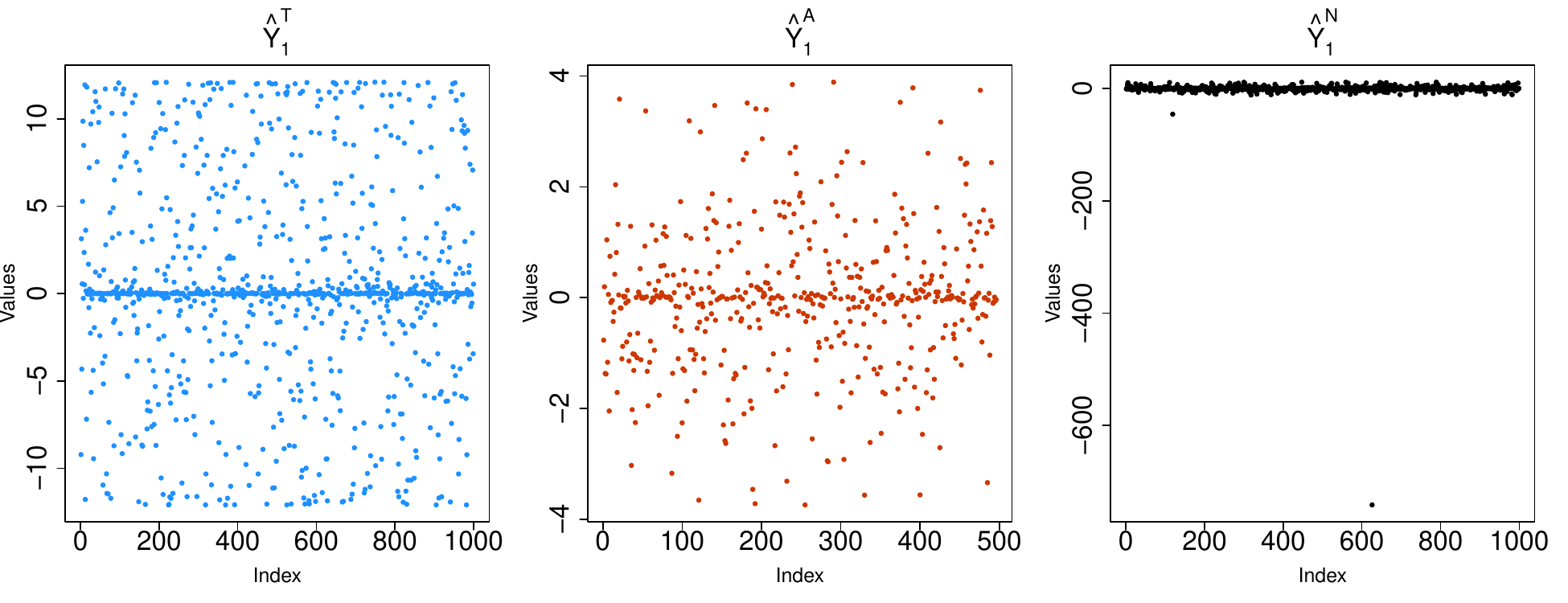}
    \caption{ Plots for $D_3$. }
    \label{D_3}
\end{figure}

\textbf{Results for $D_3$}:
It is evident from Table \ref{Table:D3} and Figure \ref{D_3} that Algorithm~\ref{Algorithm} provides good estimates of true values around 0. However, the values at both extremes of the true values are not estimated well. Here also, the values in $\hat{Y}_1^N$ are too noisy.        

\begin{table}[htbp!]
\centering
\caption{Summary statistics for $D_3$}
\label{Table:D3}
\begin{tabular}{lcccccc}
\toprule
\textbf{Values} & \textbf{Min.} & $\boldsymbol{Q_1}$ & \textbf{Median} & \textbf{Mean} & $\boldsymbol{Q_3}$ & \textbf{Max.} \\
\midrule
$\hat{Y}_1^T$    & -12.10 & -1.602 & 0.0007 & 0.2759 & 2.305 & 12.10 \\
$\hat{Y}_1^A$    & -3.7409 & -0.5490 & -0.0002 & -0.0060 & 0.5035 & 3.8894 \\
$\hat{Y}_1^N$    & -741.4168 & -0.9544 & -0.0012 & -0.8033 & 0.7530 & 11.2697 \\
\bottomrule
\end{tabular}
\end{table}

\section*{Acknowledgments}
The authors thank the editor Professor Arnak Dalalyan, an anonymous associate editor, and two anonymous reviewers for their excellent comments and constructive criticisms, which improved the article significantly.
Subhra Sankar Dhar gratefully acknowledges his core research grant (Grant Number: CRG/2022/001489), from the Government of India.

\vspace{0.25in}

\noindent {\bf Code availability:} The codes are available to the first author upon request.

\bibliographystyle{apalike} 
\bibliography{main}       
\end{document}